\documentclass[a4paper]{amsart}
\usepackage{amssymb}
\usepackage{stmaryrd}
\usepackage[all]{xy}
\usepackage{epsf}
\usepackage{amscd}

\newcommand{\printname}[1] {}

\newtheorem{theorem}{Theorem}
\newtheorem{main-theorem}{Main theorem}
\newtheorem{lemma}{Lemma}[section]
\newtheorem{proposition}[lemma]{Proposition}
\newtheorem{defi}[lemma]{Definition}
\newtheorem{corollary}[lemma]{Corollary}
\newtheorem{num}{}[section]
\newtheorem{ex}{Example}

\newcommand{\R}{\ensuremath{\mathbb R}}  

\newcommand{\Z}{\ensuremath{\mathbb Z}}  

\theoremstyle{remark}
\newtheorem{remark}{Remark}
\newcommand{\G}{\mathcal{G}}

\newcommand{\rmap}{\longrightarrow}

\begin{document}
\title[Integrability of Jacobi and Poisson structures]{Integrability of Jacobi structures}
\author{Marius Crainic}
\address{Depart. of Math., Utrecht University, 3508 TA Utrecht, The
  Netherlands}
\email{crainic@math.uu.nl}
\author{Chenchang Zhu}
\address{Department of Mathematics, University of California,
  Berkeley, CA 94720, U.S.A.}
\email{zcc@math.berkeley.edu}
\thanks{The first author was supported in part by a KNAW and a Miller
Research Fellowship. }

\date{\today}

\begin{abstract} We discuss the integrability of Jacobi manifolds by contact groupoids,
and then look at what the Jacobi point of view brings new into
Poisson geometry. In particular, using contact groupoids, we prove
a Kostant-type theorem on the prequantization of symplectic
groupoids, which answers a question posed by Weinstein and Xu
\cite{prequan}. The methods used are those of Crainic-Fernandes on
$A$-paths and monodromy group(oid)s of algebroids. In particular,
most of the results we obtain are valid also in the non-integrable
case.
\end{abstract}
\maketitle

\tableofcontents

\section{Introduction}

The ``integrability'' in the title refers to the global geometric structures
behind infinitesimal
data. Examples of ``integrations" are given by Lie's third theorem (which
integrates finite dimensional Lie algebras),
by Palais' work on integrability of infinitesimal Lie algebra actions \cite{Palais}, by
Weinstein's symplectic groupoids
which integrate Poisson structures \cite{alan-symp} (and variations, e.g. Dirac structures \cite{bcwz})
or by the integrability of general ``Lie brackets of geometric type'', i.e.
Lie algebroids \cite{marius}.

The structures that we want to integrate here are Lichnerowicz's {\it Jacobi
manifolds}
\cite{jacobi1}, known also as Kirillov's ``local Lie algebra structures'' on $C^{\infty}(M)$ \cite{kirillov}.

\begin{defi}
A Jacobi manifold $(M, \Lambda, R)$ is a manifold $M$ together with a
bivector field
$\Lambda$ and a vector field $R$ satisfying
\begin{equation} \label{jacobi}
[\Lambda,\Lambda]=2R\wedge\Lambda,\;\; \; \;[\Lambda,R]=0,
\end{equation}
The vector field $R$ is called the Reeb vector field of $M$.
\end{defi}

In the equation above, $[\cdot, \cdot]$ stands for the
Schouten-Nijenhuis bracket on multivector fields. Using Kirillov's
terminology, $C^{\infty}(M)$ together with the bracket $\{f, g\}=
\Lambda(df, dg)+ R(f)g- fR(g)$ is a local Lie algebra (and
$\Lambda$ and $R$ can be recovered from this bracket).

There are three types of ``extreme examples'':

1. {\it Contact manifolds:} To give a contact form $\theta$ on an
$(2n+1)$-dimensional manifold $M$, i.e. a 1-form with the property
that $\theta\wedge (d\theta)^n\neq 0$,  is equivalent to giving a
Jacobi structure with the property that $\Lambda^n\wedge R\neq 0$;
the defining formula is $i_{\theta}(\Lambda)= 0$, $i_{R}(\theta)=
1$.

2. {\it Vector fields:} Clearly, vector fields on $M$ can be seen
as Jacobi structures with vanishing bivector.

3. {\it Poisson manifolds:} Also,
a Poisson structure on $M$ is the same thing as a Jacobi structure
with vanishing Reeb vector field; then $\Lambda$ is called a Poisson
bivector. Note that multiplying a Poisson structure $\Lambda$ on $M$
by a smooth function $f$, the new structure $\Lambda_f= f\Lambda$ will no
longer be Poisson unless the Hamiltonian vector field $X_f$
(obtained by contracting $df$ and $\Lambda$) is zero;
instead, $\Lambda_f$ together with $X_f$ always defines a Jacobi structure.

There is yet another connection between
Jacobi and Poisson manifolds: to any Jacobi manifold $(M, \Lambda, R)$
one can associate a Poisson manifold \cite{jacobi1}.

\begin{defi}
\label{poissonization-form}
The poissonization of a Jacobi manifold $(M, \Lambda, R)$ is the Poisson
manifold
$M\times\mathbb{R}$ with the bivector:
\begin{equation}
\label{poissonization}\
\tilde{\Lambda}= e^{-s} \left(\Lambda+\frac{\partial}{\partial s}\wedge R
\right),
\end{equation}
where $s$ is the coordinate on $\mathbb{R}$. When $M$ is contact,
$M\times\R$ will be
called the symplectification of $M$ (since $\tilde{\Lambda}$ is
non-degenerate, it
defines a symplectic form).
\end{defi}

To understand the global picture behind Jacobi structures it is
useful to first recall what happens in the Poisson case, when one
discovers Weinstein's symplectic groupoids \cite{alan-symp}. First
of all, given a Poisson manifold $P$, one has an associated
topological groupoid $\Sigma_s(P)$ which shows up as the phase
space of the Poisson-sigma model \cite{cafe}, or as the
``cotangent monodromy groupoid'' of the Poisson manifold
\cite{marius, marius2}. We will call it {\it the symplectic
monodromy groupoid of $P$}. The terminology, and the subscript
``s'' come from the fact that, when $\Sigma_s(P)$ is smooth, then
it is naturally a {\it symplectic groupoid} of $P$, i.e. it comes
endowed with a symplectic form which is
multiplicative\footnote{recall that a form $\omega$ on a Lie
groupoid $\Sigma$ is called {\it multiplicative} if $m^*\omega=
pr_{1}^{*}\omega+ pr_{2}^{*}\omega$, where $pr_{1}, pr_{2}$ are
the projections, and $m$ is the multiplication, all defined on the
space $\Sigma_2$ of pairs of composable arrows of $\Sigma$.} and
it induces the Poisson structure on $P$.

This correspondence between Poisson manifolds and symplectic
groupoids, which is ``almost one-to-one'', is best explained
through the infinitesimal version of Lie groupoids, i.e. Lie
algebroids. At this point, let us fix some notations and basic
definitions (for more details and proper references, please see
\cite{marius}). For a groupoid $\Sigma$ over $M$ (hence $\Sigma$
is the space of arrows and $M$ is the space of objects, also
identified with the subspace of $\Sigma$ consisting of the
identity arrows $1_x$ at points $x\in M$), we denote by $\alpha$
the source map, by $\beta$ the target map, and by $m(g, h)= gh$
the multiplication (defined when $\alpha(g)=\beta(h)$). Hence
$\alpha, \beta: \Sigma\rmap M$, $m:\Sigma_2\rmap \Sigma$, where
$\Sigma_2$ is the space of composable arrows. Lie groupoids will
have smooth structure maps, $\alpha$ and $\beta$ will be
submersions (so that $\Sigma_2$ is a manifold), and, although the
base manifold $M$ and the $\alpha$-fibers $\alpha^{-1}(x)$ are
assumed to be Hausdorff, $\Sigma$ may be a non-Hausdorff manifold
(important examples come from bundles of Lie algebras and
foliations). Recall also that a Lie algebroid $A$ over $M$ is a
vector bundle together with a Lie bracket $[\cdot, \cdot]$ on the
space of sections $\Gamma(A)$ and a bundle map $\rho:
A\longrightarrow TM$ so that the Leibniz identity holds:
\[ [\xi_1, f\xi_2]= f[\xi_1, \xi_2]+ \mathcal{L}_{\rho(\xi_1)}(f)\xi_2 \]
for all sections $\xi_1$ and $\xi_2$ and all
functions $f$; here $\mathcal{L}$ stands for Lie derivatives along vector fields.
$\rho$ is called {\it the anchor} of $A$.
As in the case of Lie groups and Lie algebras, any Lie groupoid has an
associated Lie algebroid
(obtained by taking the tangent spaces along the $\alpha$-fibers
at each identity element $1_x$). However, not all Lie algebroids arise in
this way;
one says that $A$ is {\it integrable} if it comes from a Lie groupoid.
Nevertheless, any Lie algebroid $A$ has an associated {\it monodromy
groupoid} $G(A)$ (called also the Weinstein groupoid of $A$), made out of
homotopy classes of paths in the
algebroid world, and it is a topological groupoid which is the universal
candidate for integrating $A$
\cite{marius} (this will be recalled in section \ref{section3}).

Going back to Poisson manifolds, the cotangent bundle $T^*P$
of any Poisson manifold $P$
is naturally a Lie algebroid: the anchor is the contraction by
the Poisson tensor:
\begin{equation}
\label{Lambda-sharp}
\Lambda^{\sharp}: T^*M\rmap TM,\ \Lambda(\omega, \eta)= \eta(\Lambda^{\sharp}(\omega)),
\end{equation}
and the bracket is
\begin{equation}
\label{bracket-Lambda}
[\omega, \eta]_{\Lambda}= \mathcal{L}_{\Lambda^{\sharp}\omega}\eta- \mathcal{L}_{\Lambda^{\sharp}\eta}\omega- d\Lambda (\omega, \eta).
\end{equation}
Note that the bracket is uniquely determined by the Leibniz identity and
$[df, dg]=d\{f, g\}$. With these, the symplectic monodromy groupoid of $P$ is
just
\[ \Sigma_s(P)= G(T^*P), \]
the monodromy groupoid associated to $T^*P$.

For Jacobi manifolds, there are partial results which are similar to those
in the Poisson case,
and the necessary objects are already known \cite{ dazord, KSB}.

\begin{defi}
A contact groupoid over a
manifold $M$ is a Lie groupoid $\Sigma$ over $M$ together with a contact
form $\theta$, and a function
$r$ on $\Sigma$, with the property that $\theta$ is $r$-multiplicative in
the sense that
\begin{equation}
\label{contact gpd}
m^* \theta = pr_2^*(e^{-r}) \cdot pr_1^* \theta + pr^*_2 \theta .
\end{equation}
(where, as in the previous footnote, $m$ is the multiplication and $pr_i$
are the projections).
The function $r$ is called the Reeb function, or the Reeb cocycle of
$\Sigma$.
\end{defi}

For a discussion on the non-symmetry of the previous equation, and versions which
use the point of view of contact hyperplanes instead of contact forms, please see
Example \ref{loclly-conf} in the last section.

As we shall explain, Reeb cocycles come from integrating Reeb
vector fields of Jacobi manifolds. The term ``cocycle'' comes from
the fact that the definition above forces the cocycle condition
$r(gh)= r(g)+ r(h)$, whenever $gh$ is defined. This implies that
the base $M$ of a contact groupoid has an induced Jacobi structure
\cite{KSB} (and it will also follow from the next section).

As in the Poisson case, any Jacobi manifold has an associated Lie algebroid \cite{KSB},
hence a monodromy groupoid.

\begin{defi}
\label{def-alg-jacobi}
The Lie algebroid of the Jacobi manifold $(M, \Lambda, R)$
is $T^*M\oplus \mathbb{R}$, with the anchor
$\rho: T^*M\oplus\mathbb{R}\longrightarrow TM$ given by
\[ \rho(\omega, \lambda)= \Lambda^{\sharp}(\omega)+ \lambda R ,\]
and the bracket
\begin{align}
& [(\omega, 0), (\eta, 0)]= ([\omega, \eta]_{\Lambda}, 0)- (i_R(\omega\wedge\eta), \Lambda(\omega, \eta))\nonumber\\
& [(0, 1), (\omega, 0)]= (L_{R}(\omega), 0), \nonumber
\end{align}
(where $\Lambda^{\sharp}$ and $[\cdot, \cdot]_{\Lambda}$ are given by (\ref{Lambda-sharp}) and (\ref{bracket-Lambda}) above).
The associated groupoid
\[ \Sigma_{c}(M)= G(T^*M\oplus \mathbb{R}) ,\]
is called the contact monodromy groupoid of the Jacobi manifold $M$.
We say that $M$ is integrable as a Jacobi manifold if the associated algebroid
$T^*M\oplus \mathbb{R}$ is integrable (or, equivalently, if $\Sigma_{c}(M)$ is smooth, cf. \cite{marius}).
\end{defi}

The fact that $R$ can be integrated to define a multiplicative
function on $\Sigma_{c}(M)$ (which will be explained in detail in
the main body), together with the local result of Dazord
\cite{dazord}, suggests that $\Sigma_{c}(M)$ is a contact groupoid
whenever it is smooth. Our first main result proves that this is
indeed the case, and also describes the relation between the
integrability of $M$ and of its Poissonization. To describe this
relation at the groupoid level, one remarks that any
multiplicative function $r$ on a groupoid $\Sigma$ over $M$ can be
used to define an ``$r$-twisted multiplication by $\R$'', which is
a groupoid $\Sigma\times_r \R$ over $M\times\R$ (cf. Definition
\ref{multiply-reals}). We then have:

\begin{main-theorem}
\label{main}
For any Jacobi manifold $M$,
\begin{enumerate}
\item[(i)] there is an isomorphism of topological groupoids
\[ \Sigma_{s}(M\times\R)\cong \Sigma_{c}(M)\times_{r}\R ,\]
and $M$ is integrable if and only if  the Poisson manifold $M\times \R$ is
integrable.
\item[(ii)] $M$ is integrable if and only if $\Sigma_{c}(M)$
is smooth. Moreover, in this case $\Sigma_{c}(M)$ has a natural structure of
contact groupoid.
\end{enumerate}
\end{main-theorem}

Next, we concentrate on Poisson geometry by viewing Poisson manifolds as
Jacobi
ones with vanishing Reeb vector fields. In other words, we look at what the
``Jacobi point of view''
(rather than Jacobi structures) brings new into Poisson geometry. First of
all, it shows that a Poisson manifold $P$ comes not only with the symplectic monodromy
groupoid $\Sigma_s(P)$, but also with a contact monodromy groupoid $\Sigma_{c}(P)$.
Of course, it is not a surprise that
the relation between the two heavily depends on the Poisson geometry of $P$.
However, this relation also provides new insights, e.g. into the geometric
prequantization of Poisson manifolds.
First of all, we concentrate on describing the relation between the two
groupoids. Here we restrict to the case where $P$ is integrable as a
Poisson manifold (the general case is treated in Section \ref{Poisson-I}). We have a
bundle
of groups over $P$, $\mathcal{P}_{\Lambda}$, whose fiber at $x\in P$
is the period group of the restriction of $\Omega$ to the $\alpha$-fiber
${\alpha^{-1}(x)}$,
where $\Omega$ is the symplectic form of $\Sigma_s(P)$.
We also define $\mathcal{G}_{\Lambda}$ as the quotient of the trivial bundle
with
fiber $\R$ by $\mathcal{P}_{\Lambda}$.

\begin{main-theorem}
\label{main2}
For any Poisson manifold $P$, there is a short exact sequence of
topological
groupoids
\[ 1\rmap \mathcal{G}_{\Lambda}\rmap \Sigma_{c}(P)\rmap \Sigma_{s}(P)\rmap 1 .\]
Moreover, if $P$ is integrable as a Poisson manifold, the following are equivalent
\begin{enumerate}
\item[(i)] $P$ is integrable as a Jacobi manifold.
\item[(ii)] $\mathcal{P}_{\Lambda}$ is smooth.
\item[(iii)] $\mathcal{G}_{\Lambda}$ is smooth.
\end{enumerate}
\end{main-theorem}

Next, we will use the contact groupoid $\Sigma_{c}(P)$ to clarify
the prequantization
of the symplectic groupoid $\Sigma_s(P)$. In particular, we will see that
any prequantization
is a contact groupoid. As a simplified theorem that we can state in this
introduction, we mention here:

\begin{main-theorem}
\label{main3}
For any integrable Poisson manifold $P$, the following are
equivalent:
\begin{enumerate}
\item[(i)] $\Sigma_{s}(P)$ is prequantizable.
\item[(ii)]  $\mathcal{P}_{\Lambda}\subset P\times \mathbb{Z}$.
\end{enumerate}
Moreover, if $\Sigma_s(P)$ is Hausdorff, the conditions above are also
equivalent to 
\begin{enumerate}
\item[(iii)] The symplectic form of $\Sigma_s(P)$ is integral.
\end{enumerate}
In this case, the prequantization $\tilde{\Sigma}$ is Hausdorff.
\end{main-theorem}

And, finally, we will describe the connection with the Van Est map. More precisely, the Poisson tensor
can be viewed as an algebroid 2-cocycle on $T^*M$, and it makes sense to ask when it is integrable (i.e.
when it comes from a 2-cocycle on the symplectic groupoid, via the Van Est map). We will show
(compare with the previous theorem):

\begin{main-theorem}
\label{main4} Let $(P, \Lambda)$ be an integrable Poisson
manifold, and $\Sigma_s(P)$ is Hausdorff. Then the following are
equivalent:
\begin{enumerate}
\item[(i)] $\Lambda$ is integrable as an algebroid cocycle.
\item[(ii)] The symplectic form of $\Sigma_s(P)$ is exact.
\item[(iii)] $\mathcal{P}_{\Lambda}= 0$.
\end{enumerate}
Moreover, in this case $\Sigma_c(P)$ is isomorphic to the product $\Sigma_s(P)$
with $\mathbb{R}$, with the multiplication twisted by a cocycle on $\Sigma_s(P)$
integrating $\Lambda$.
\end{main-theorem}

The paper is organized as follows. In the first section we give more details on the poissonization process,
including a groupoid version. Next, there is one section devoted to each of the main theorems, which provides
the details on the definitions, more precise statements, and the proofs. In the last section we look at several particular cases
and examples.

{\it Late comment:} Finally, we would like to mention that, at the time the first version of this paper was written,
we were not aware of several related works. Most notably, the non-trivial remark that the integrability
of Jacobi structures is intimately related to prequantization already appears in \cite{dazord2}.
Also, the idea of viewing the Reeb vector field of a Jacobi manifold as
an algebroid cocycle already appears in \cite{igma} (where it plays a central role). However, regarding the
main results, there is hardly any overlap, and our methods (and point of view) is very different
from the existing ones. One of the reasons is due to the fact that we use the $A$-path approach
of \cite{marius} (and this has been proven to be more powerful also in the Poisson case \cite{marius2}).

{\bf Acknowledgments:} The second author would like to thank her
advisor Alan Weinstein for very helpful comments and discussions.

\section{Poissonization and homogeneity}
\label{Poissonization and homogeneity}

\begin{num}\underline{{\bf The manifold case \cite{jacobi1, jacobi2}.}}\rm \
\label{hom-Ps}
A {\it homogeneous Poisson manifold}
is a Poisson manifold $(P, \tilde{\Lambda})$ together with a vector
field $Z$, called the {\it homogeneous vector field} of $P$, with the
property that
\[ \tilde{\Lambda}=-  \mathcal{L}_{Z}(\tilde{\Lambda}).\]
If $\tilde{\Lambda}$ comes from a
symplectic form $\omega$, then the equation above becomes $\omega=
\mathcal{L}_{Z}(\omega)$,
and one calls $(P, \omega, Z)$ a {\it homogeneous symplectic manifold}.
Recall also that the {\it poissonization} of
a Jacobi manifold $(M, \Lambda, R)$ is the Poisson manifold
$(M\times \mathbb{R}, \tilde{\Lambda})$, where $\tilde{\Lambda}$
is given by the formula (\ref{poissonization}) in the introduction.
We have (cf. subsection 11 in \cite{jacobi1}, or the theorem on pp. 443 of \cite{jacobi2}):
\end{num}

\begin{proposition}
\label{prop1-manifold}
The poissonization construction defines a one-to-one correspondence
between Jacobi structures on $M$ and homogeneous Poisson structures on
$M\times \R$ with homogeneous vector field $\frac{\partial}{\partial s}$.

Moreover, when restricted to contact manifolds, this induces a one-to-one
correspondence between contact structures on $M$ and homogeneous symplectic
structures on $M\times \R$ with homogeneous vector field
$\frac{\partial}{\partial s}$.
\end{proposition}

\begin{proof} The proof is based on some simple remarks which are
interesting on their own.
First, for the inverse of the poissonization
construction we remark that, given a homogeneous Poisson structures $\tilde{\Lambda}$ on
$M\times \R$ with homogeneous vector field $\frac{\partial}{\partial s}$,
one has an induced Jacobi structure on $M$ with
\[ \Lambda= (pr_{1})_{*} (e^s \tilde{\Lambda}) , \;\;\;\;\; R= e^s\tilde{\Lambda}^{\sharp} (d s), \]
where $pr_1: M\times \R \longrightarrow M$ is the projection.
Next, when the Jacobi manifold $M$ is actually contact with contact form
$\theta$,
then its poissonization is actually symplectic, with the symplectic form
\[ \omega = d(e^s pr_{1}^{*} \theta).\]
And, finally, if
$(M\times\mathbb{R}, \omega)$ is symplectic and
$\mathcal{L}_{\frac{\partial}{\partial s}} \omega = \omega$,
then the induced Jacobi structure on $M$ comes from the contact form
$\theta =  \iota (\frac{\partial}{\partial s}) \omega$.
\end{proof}

\begin{num}\underline{{\bf The groupoid case.}}\rm \
We now discuss the groupoid version of the proposition above.
Corresponding to Poisson manifolds are symplectic groupoids,
i.e. Lie groupoids $\Sigma \underset{ \beta}
{\overset{\alpha} {\rightrightarrows} }P$ equipped with a
symplectic form $\omega$ on $\Sigma$ which is multiplicative,
i.e. which satisfies
\begin{equation} \label{symp gpd}
m^* \omega = pr_1^* \omega + pr^*_2 \omega,
\end{equation}
where the equation is on the space $\Sigma_2$ of pairs of
composable elements, $pr_j$ is the projection into the $j$-th
component, and $m$ is the multiplication. Recall also that a
multiplicative vector field on $\Sigma$ consists of a vector field
$Z$ on $\Sigma$ and a vector field $Z_0$ on the base manifold $M$,
with the property that the flow $\phi_{Z}^{t}$ is a local Lie
groupoid morphism over the flow $\phi_{Z_0}^{t}$. That means that
for any arrow $g$ from $x$ to $y$ so that $g'= \phi_{Z}^{t}(g)$ is
defined, $x'= \phi_{Z_0}^{t}(x)$ and $y'= \phi_{Z_0}^{t}(y)$ are
defined too and $g'$ is an arrow from $x'$ to $y'$, and the
multiplicativity condition $\phi_{Z}^{t}(gh)= \phi_{Z}^{t}(g)
\phi_{Z}^{t}(h)$ holds whenever the right hand side is defined.
These conditions can be reformulating by saying that $Z:
\Sigma\rmap T\Sigma$ is a groupoid morphism with base map $Z_0:
M\rmap TM$, where $T\Sigma$ is with the induced structure of
groupoid over $TM$ (for details, see \cite{m-m}. Since $Z_0$ can
be recovered from $Z$ (push down along the source map), we simply
say that $Z$ is a multiplicative vector field.

Finally,
\end{num}

\begin{defi}
A homogeneous symplectic groupoid $(\Sigma, \omega, Z)$ is a
symplectic groupoid together with a multiplicative vector field $Z$,
with the property that $\mathcal{L}_{Z}(\omega)= \omega$.
\end{defi}

Next, we need the groupoid version of ``multiplying with the reals'' that
appears in the poissonization procedure above.

\begin{defi}
\label{multiply-reals}
Given a groupoid $\Sigma$ (Lie or not), and a multiplicative function $r$ on
$\Sigma$,
we define the groupoid $\Sigma\times_{r} \R= \Sigma\times \R$ over
$M\times\R$,
with source $\alpha$, target $\beta$ and multiplication given by
\begin{equation}
\label{extension}
\alpha(g, s)= (\alpha(g), s), \beta(g, s)= (\beta(g), s- r(g)), (g_1, s_1)
(g_2, s_2)= (g_1g_2, s_2).
\end{equation}
\end{defi}

Actually, it is easy to check that $\Sigma\times_{r} \R$ is a groupoid if
and only if $r(g_1g_2)= r(g_1)+r(g_2)$, i.e.
$r$ is multiplicative,

\begin{proposition}\label{prop1-groupoid}
Let $\Sigma$ be a Lie groupoid endowed with a smooth multiplicative function
$r$. Then there is a one-to-one  correspondence between
contact groupoid structures on the Lie groupoid $\Sigma$ with Reeb function
$r$,
and homogeneous symplectic groupoid structures on
the Lie groupoid $\Sigma\times_r\R$ with homogeneous vector field
$\frac{\partial}{\partial s}$.
\end{proposition}

\begin{proof} By the second part of the previous proposition we are left
with showing that,
if the groupoid $\Sigma$ comes equipped with the multiplicative function $r$
and a contact form $\theta$,
and if $\omega$ is the induced symplectic structure on
$\Sigma\times\mathbb{R}$, then
the multiplicativity of $\omega$ (i.e. equation (\ref{symp gpd})) is
equivalent
to the $r$-multiplicativity of $\theta$ (i.e. equation (\ref{contact gpd})).
For this one recalls that $\omega=d(e^s\theta)$, and one remarks that
the space of composable pairs of arrows in $\Sigma\times_{r}\R$ can be
identified with
$\Sigma_2\times\R$ by
\[ \left( (g, s-r(g')),(g', s)\right)\mapsto ( g, g', s).\]
Taking the component of \eqref{symp gpd} containing $ds$, we obtain
\eqref{contact gpd} immediately.
The other direction follows by multiplying by $e^{s}$ and taking
derivatives.
\end{proof}

\begin{defi}
\label{def-symplectification}
Given a contact groupoid $\Sigma$ with Reeb function $r$, the associated
symplectic
groupoid $\Sigma\times_{r}\R$ is called {\it the symplectification} of
$\Sigma$.
\end{defi}

\begin{num}\underline{{\bf Compatibility.}}\rm \
We now point out the compatibility between the correspondences
described by the previous two propositions. Recall that, given a
symplectic groupoid $\Sigma$ over $P$, there is an induced Poisson
structure on $P$, uniquely determined by the property that the
source map is a Poisson map\footnote{This construction actually
gives one-to-one correspondences between symplectic groupoids over
$P$ and integrable Poisson structures on $P$.} (i.e. preserves the
Poisson bivector). A similar result holds for contact groupoids
and Jacobi manifolds: 
\end{num}

\begin{lemma}
\label{induced-on-base} Given a contact groupoid $\Sigma$ over
$M$, there exists an unique Jacobi structure with the property
that the source map $\alpha: \Sigma\rmap M$ is a Jacobi map. In
this case, we call $\Sigma$ a contact groupoid of the Jacobi
manifold $M$.
\end{lemma}

The proof of this lemma will also show the following which,
although straightforward, proves that the correspondence of
Proposition \ref{prop1-groupoid} implies that of Proposition
\ref{prop1-manifold}.

\begin{proposition}
\label{identify}
If $(\Sigma, \theta, r)$ is a contact groupoid
over $M$, and $(\Sigma\times_{r}\mathbb{R}, \omega)$ is the associated
homogeneous symplectic groupoid, then the poissonization of the Jacobi
structure
induced on $M$ (by the contact groupoid $\Sigma$) coincides with the
Poisson structure induced on $M\times\mathbb{R}$
(by the symplectic groupoid $\Sigma\times_{r}\mathbb{R}$).
\end{proposition}

\begin{proof} (of Lemma \ref{induced-on-base} and of Proposition
\ref{identify})
The uniqueness in the lemma is clear since $\alpha$ is a submersion. We
prove the rest. Remark that, in general, a Poisson tensor $\tilde{\Lambda}$ is homogeneous
with respect to a vector field $Z$ if and only if $\phi^{t}_{Z}$ maps $\tilde{\Lambda}$ into
$e^t\tilde{\Lambda}$. On the other hand, if $\Gamma$ is a homogeneous symplectic groupoid over $P$ with
homogeneous vector field $Z$, we know that the flow of
$Z$ is a (local) groupoid homomorphism over the flow of $Z_0$. We deduce that
the induced Poisson structure on $P$ is a homogeneous one, with homogeneous vector field $Z_0$.
Now, given a contact groupoid $\Sigma$ over $M$, we form the homogeneous symplectic
groupoid $\Sigma\times_{r}\R$, and it follows that the induced Poisson structure on
$M\times \R$ is homogeneous with vector field $\frac{d}{ds}$. Hence it comes
from a Jacobi structure on $M$. One still needs to remark that, by the
correspondence of Proposition \ref{prop1-manifold}, Jacobi maps correspond to Poisson maps
(and we apply this to the source map).
\end{proof}

\section{Symplectic and contact monodromy groupoids}
\label{section3}

\begin{num}\underline{{\bf The main theorem of the section.}}\rm \
In this section we investigate the effect that
the poissonization process
has on the monodromy groupoids. We first recall the construction of the
monodromy groupoid $G(A)$ associated to a Lie algebroid $A$. As mentioned
in the introduction,
when applied to the algebroid $T^*P$ of a Poisson manifold $P$ and to the
algebroid $T^*M\oplus\mathbb{R}$ of a Jacobi manifold $M$, one defines
\[ \Sigma_s(P)= G(T^*P), \quad \Sigma_c(M)= G(T^*M\oplus\mathbb{R}), \]
called the symplectic monodromy groupoid of $P$, and the contact monodromy
groupoid of $M$, respectively. The purpose of this section will then be to
prove the following stronger version of the Main Theorem \ref{main} stated in
the introduction.
\end{num}

\begin{theorem} \label{basic correspondence}
Let $M$ be a Jacobi manifold with Reeb vector field $R$, let $\Sigma_{c}(M)$
be the contact monodromy groupoid
of $M$ and let $\Sigma_s(M\times\R)$ be the symplectic monodromy groupoid of
the Poissonization of $M$. Then:
\begin{enumerate}
\item[(i)] By integration, the Reeb vector field $R$ induces a
multiplicative
function $r$ on $\Sigma_{c}(M)$.
\item[(ii)] One has an isomorphism of topological groupoids
\begin{equation}
\label{groupoid-level}
\Sigma_{s}(M\times\mathbb{R})\cong \Sigma_{c}(M)\times_{r} \mathbb{R} ,
\end{equation}
(where $r$ is the one from (i), and ``$\times_{r}$'' was introduced in
Definition \ref{multiply-reals}).
\item[(iii)] $M$ is integrable as a Jacobi manifold if and only if
$M\times\R$ is integrable
as a Poisson manifold.
\item[(iv)] In the integrable case, $\Sigma_{c}(M)$ is the source-simply connected contact groupoid of $M$ with Reeb function $r$, and its symplectification (cf. Definition \ref{def-symplectification})
is isomorphic to the symplectic groupoid $\Sigma_{s}(M\times \R)$
of Poisson manifold $M\times \R$.
\end{enumerate}
\end{theorem}

\begin{num}\underline{{\bf Monodromy groupoids.}}\rm \
We now recall from \cite{marius} the construction of the monodromy groupoid
$G(A)$ associated to a Lie algebroid $A$.
\end{num}
\begin{defi} \label{apath}
Given a Lie algebroid $A \overset{\pi} {\longrightarrow} M$ with
anchor $\rho : A \longrightarrow TM$, an {\em $A$-path} of $A$ is
a $C^1$ map $a$: $I=[0, 1]
\longrightarrow A$ with the property that
\[ \rho (a(t)) = \frac{d\gamma }{d t} (t),\]
where $\gamma(t)=\pi \circ a(t)$ is called the base path of
$a$.
We denote by $P_a(A)$ the set of all $A$-paths of $A$.
\end{defi}

The choice of the order of smoothness is not very important, and
we choose it finite in order to work with Banach manifolds
and not with Frechet ones. In particular, the space
$P(A)$ of all paths in $A$ will be a Banach manifold modelled by the Banach
space
$P(\R^n)=C^1(I, \R^n)$, and $P_a(A)$ will be a Banach submanifold.
For more details, see \cite{marius}.

\begin{defi}
\label{homotopy-def}
\cite{marius} Let $\nabla$ be a connection
on $A$
with torsion $T_{\nabla}$ defined as
\[ T_{\nabla} ( \alpha , \beta) = \nabla_{\rho(\beta)} \alpha -
\nabla_{\rho(\alpha)} \beta + [\alpha, \beta], \]
and let $\partial_{t}$ be the induced derivative operator (which associates
to a path $a=a(t)$ in $A$,
the path in $A$ $\partial_{t}a$ which is the $\nabla$-horizontal component
of $\frac{da}{dt}$).
An $A$-homotopy is a family $a_{\epsilon}(t)=a(\epsilon, t)$ of $A$-paths of
class $C^2$ in $\epsilon$
with the property that their base paths
$\gamma_{\epsilon}(t)=\gamma(\epsilon,t)$ have fixed end points,
and the solution $b= b(\epsilon, t)$ of the equation
\begin{equation}
\label{homotopy}
\partial_t b -\partial_{\epsilon }a = T_{\nabla} (a, b),
\;\;\;\;\; b(\epsilon ,0)=0
\end{equation}
satisfies $b(\epsilon,1)=0$ for all $\epsilon$. We say that two $A$-paths
$a_0$ and $a_1$ are homotopic, and write $a_0\sim a_1$, if there exists an
$A$-homotopy $a(\epsilon, t)$ with
the property that $a_{i}(t)= a(i, t)$, $i= 0,1$.
\end{defi}

Recall \cite{marius} that the solution $b$ of the previous equation (hence
also the homotopy relation)
does not depend on the choice of the connection $\nabla$. Intuitively,
$A$-homotopies are ``algebroid homotopies with fixed end-points'', and the
equation (\ref{homotopy})
above is just the algebroid version of the equation (in $\R^n$)
$\frac{d}{dt}\frac{d}{d\epsilon}= \frac{d}{d\epsilon}\frac{d}{dt}$
which can be used to compute $b= \frac{d\gamma}{d\epsilon}$ from $a=
\frac{d\gamma}{dt}$.

Finally, the monodromy groupoid of $A$ is defined as
\[ G(A):= \big( P_a(A)/\sim \big)
\underset{\beta}{\overset{\alpha}{\rightrightarrows}} M.\]
The source and target maps $\alpha$ and $ \beta$ are given by
\[ \alpha([a]) = \gamma(0), \;\;\;\;\; \beta([a]) = \gamma(1), \]
where $\gamma$ is the base path of $a$.
The multiplication is basically the concatenation of paths:
\[
a\odot b(t)\equiv \left\{
\begin{array}{ll}
2b(2t),\qquad& 0\le t\le \frac{1}{2}\\ \\
2a(2t-1),\qquad & \frac{1}{2}< t\le 1
\end{array}
\right.
\]
Strictly speaking, this forces us to consider the pathwise smooth
$A$-paths. Instead, since reparametrization does not affect the
homotopy class \cite{marius}, we can first reparametrize $A$-paths
by a cut-off function $\tau \in C^{\infty} (I, I)$ (the
reparametrization of $a$ is $a^{\tau}(t)=\tau'(t)a(\tau(t))$, and
then define $[a]\cdot [b]= [a^{\tau}\odot b^{\tau}]$.

Since $G(A)$ is the quotient of the Banach manifold $P_a(A)$, it follows
that
$G(A)$ is a topological groupoid, and we can talk unambiguously about its
smoothness:
we are only interested on smooth structures for which the projection from
$P_{a}(A)$ onto $G(A)$
is a submersion, and there is at most one such structure. It then follows that
$A$ is integrable if and only
if $G(A)$ is smooth, in which case $G(A)$ will be the unique Lie groupoid
integrating $A$ which
has simply connected $\alpha$-fibers \cite{marius}.

\begin{num}\underline{{\bf Passing to 1-cocycles.}}\rm \
To prove (i), (ii) and (iii) of Theorem \ref{basic correspondence}, it is
useful  to
concentrate on the algebroid
\[ A= T^*M\oplus \R, \]
to identify the Reeb vector field with the section $R\in \Gamma(A^*)$ which
vanishes in the $\R$-direction:
\[ R(\omega, \lambda)= \omega(R) ,\]
and to remark that $R$ becomes an algebroid 1-cocycle, i.e.
\[ R([\alpha, \beta])= \mathcal{L}_{\rho{\alpha}}(R(\beta))-
\mathcal{L}_{\rho{\beta}}(R(\alpha))\]
for all $\alpha, \beta\in \Gamma (A)$.
\end{num}

\begin{defi} Given a Lie algebroid $A$ and a 1-cocycle $R\in
\Gamma(A^*)$, define $A\times_{R}\R$
as the algebroid over $M\times\R$, which, as a vector bundle, is the
pull-back of $A$ to $M\times\R$
along the projection, has the anchor $\rho_{R}: A\times_{R}\R\rmap
T^{*}(M\times\R)$,
\[ \rho_{R}(\alpha)=  \rho (\alpha)-
R(\alpha)\frac{\partial}{\partial s} ,\] and the Lie bracket:
\[ [\alpha, \beta]_R= [\alpha, \beta] -
R(\alpha)\frac{\partial \beta}{\partial s} +R(\beta)\frac{\partial
\alpha}{\partial s},
\]
where $\rho$ and $[\cdot, \cdot]$  are anchor and Lie bracket on $A$.
\end{defi}

Via the remark that 1-cocycles $R\in \Gamma(A^*)$ are the same thing as Lie algebroid actions
of $A$ on $M\times \mathbb{R}$, the algebroid $A\times_{R}\R$ is the standard action algebroid
(or pull-back algebroid) associated to the action (see e.g. \cite{MaHi, m-m}).
The point of this definition is that it allows us to
include the algebroid $T^{*}(M\times \R)$ into
the general picture:

\begin{lemma}
\label{conceptual}
For the Lie algebroid $A= T^*M\oplus\R$, and the Reeb vector field
$R$ viewed as a 1-cocycle of $A$, one has an isomorphism of algebroids
\[ T^{*}(M\times\R) \cong A\times_{R} \R .\]
\end{lemma}
\begin{proof}
One uses the bundle isomorphism $\psi(v, t)=(e^{-t} v, t)$.
\end{proof}

We now have the following general result:

\begin{proposition}
\label{algebroid-ref}
Let $A$ be a Lie algebroid, and let $R\in \Gamma(A^*)$ be an
1-cocycle.
Then
\begin{enumerate}
\item[(i)] The integral
\[ r(a)= \int_a R:= \int_{0}^{1} \langle R(\gamma(t)), a(t)\rangle dt \]
only depends on the homotopy class of the $A$-path $a$, and induces
a multiplicative function $r$ on $G(A)$.
\item[(ii)] There is an isomorphism of topological groupoids
\[ G(A\times_{R}\R)\cong G(A)\times_{r} \R .\]
\item[(iii)] $A$ is integrable if and only if $A\times_{R}\R$ is. In this
case, the previous isomorphism is a Lie groupoid isomorphism.
\end{enumerate}
\end{proposition}

\begin{proof}
That $\int_{a}R$ only depends on the homotopy class of $a$ has
been proven
for Lie algebroids coming from Poisson manifolds in \cite{marius2}, and
exactly the same proof
applies in general.
That $r$ is multiplicative is clear from the additivity of integration. Next,
it is easy to see that an $A$-path
of $A\times_{R}\R$ is the same thing as an $A$-path $a$ of $A$, together with a
path $\gamma_0$ in $\R$, satisfying
\begin{equation}
\label{definingeq}
\frac{d\gamma_0}{dt}= - R(a(t)) .
\end{equation}
In turn, this determines $\gamma_0$ from the initial point $s= \gamma_0(0)$.
This defines a bijection
$P_{a}(A\times_{R}\R)\cong P_{a}(A)\times\R$ ($(a, \gamma_0)$ corresponds to
$(a, s)$), which is clearly smooth.
Moreover, choosing a connection $\nabla$ on $A$ and the pull-back
$\tilde{\nabla}$ on $A\times_{R}\R$ in order
to write the homotopy equations (\ref{homotopy}) for the two algebroids, it
is straightforward to see that
this correspondence preserves the homotopy, hence it descends to the
isomorphism of topological spaces:
\[ G(A\times_{R}\R)\cong G(A)\times \mathbb{R} .\]
It is straightforward to identify the groupoid structure on the
right hand side with $G(A)\times_{r} \R$. For instance, the source
and the target of $([a], s)$ will be $(\gamma_1(0), \gamma_0(0))=
(\alpha([a]), s)$, and $(\gamma_1(1), \gamma_0(1))= (\beta([a]),
s- r(a))$, respectively, where we have integrated
(\ref{definingeq}) to compute $\gamma_0(1)$.

For (iii) we use again $\nabla$. We can talk about geodesic $A$-paths, and
define
the exponential map
\[ exp_{\nabla}: A\rmap G(A)\]
which associates to $v\in A$ the homotopy class of the geodesic $A$-path
that starts at $v$. As
in the classical case, $exp_{\nabla}$ is defined only for small enough
$v$'s, but we
make an abuse of notation and write as if it was defined globally.
The point is that $A$ is integrable if and only if, locally around the zero
points, $exp_{\nabla}$
is injective. This has been explained in \cite{marius}, to which we refer
also for more details
on the exponential map. Then, for the first part of (iii), it suffices to
remark that,
after the identification $G(A\times_{R}\R)\cong G(A)\times_{r} \R$,
$exp_{\tilde{\nabla}}$ is identified with $exp_{\nabla}\times id$.
Also, since the smooth structure on $G(A)$ is constructed with the help of
the exponential map
(and that is why we need it to be injective), the last part of (iii)
follows.
\end{proof}

\begin{num}\underline{{\bf Proof of Theorem \ref{basic correspondence}.}}\rm \
By Lemma \ref{conceptual}, (i), (ii) and (iii) in Theorem
\ref{basic correspondence} follow from Proposition
\ref{algebroid-ref}. We are now left with part (iv). Since
$\Sigma_s(M\times\R)$ is a symplectic groupoid, the isomorphism
(\ref{groupoid-level}) makes $\Sigma_c(M)\times_r\R$ into a
symplectic groupoid with symplectic form denoted by $\omega$.
Using Proposition \ref{prop1-groupoid}, it suffices to show that
$\frac{d}{ds}$ is multiplicative and $(\Sigma_c(M)\times_r\R,
\omega)$ is homogeneous with respect to the vector field
$\frac{d}{ds}$. We only have to show that $\mathcal{L}_{\frac{d}{d
s}} \omega =\omega$. We will make some general remarks. First of
all, if $(P, \Lambda)$ is a Poisson manifold, and $\lambda$ is a
non-zero real, it is immediate from the construction of the
symplectic form $\omega_{\Lambda}$ on $\Sigma_s(P, \Lambda)$ (see
\cite{marius2}) that
\[ \Sigma_s(P, \lambda\Lambda)= \Sigma_s(P, \Lambda),\quad
\omega_{\lambda\Lambda}= \lambda^{-1} \omega_{\Lambda} .\]
Next, if $\phi: (M_1, \Lambda_1)\rmap (M_2, \Lambda_2)$ is an isomorphism
of
two Poisson manifolds, then it induces an isomorphism of algebroids
$\phi_{*}:  T^*M_1\rmap T^*M_2$ which, on each fiber, is given by the
inverse $(d\phi)^{-1}$
of the differential of $\phi$. In the integrable case, it induces a map
$\phi_{*}: \Sigma_s(M_1)\rmap\Sigma_{s}(M_2)$ of symplectic groupoids.
Applying this to $\phi_u: M\times\R\rmap M\times\R$,
$\phi_{u}(x, s)= (x, s+u)$, $\Lambda_1= \tilde{\Lambda}$, $\Lambda_2=
e^{-u}\tilde{\Lambda}$,
and using also the previous remark, we obtain
\[ (\phi_{u})_{*}: (\Sigma_{s}(M\times\R), \omega)\rmap
(\Sigma_{s}(M\times\R), e^u\omega).\] After the identification
from (iii), we see that $(\phi_{u})_{*}(g, s)= (g, s+u)$ is the
flow of $\frac{d}{ds}$. Hence $\frac{d}{ds}$ is multiplicative.
Taking derivatives in $(\phi_{u})_{*}(\omega)= e^{u}\omega$ (with
respect to $u$, at the origin), we obtain the desired equation
$\mathcal{L}_{\frac{d}{d s}} \omega =\omega$.
\end{num}

\begin{remark} \rm \
Let us go back to Theorem \ref{basic correspondence} and conclude with the
explicit formulas. First of all, the Reeb function on $\Sigma_{c}(M)$ is
given by
\[ r([a])= \int_{a} R:= \int_{0}^{1} \langle R(\gamma(t)), a_1(t)\rangle dt
,\]
for any $A$-path $a= (a_1, a_0)$ of $T^*M\oplus\R$ with base path $\gamma$.
Also, the isomorphism in (iii) comes from a diffeomorphism at the level
of $A$-paths,
\begin{equation} \label{corr}
\begin{array}{ccc}
P_a(T^*(M \times \R)) & \longrightarrow &  P_a(T^*M
\oplus_M \R) \times \R \\
\tilde{a}_1+ \tilde{a}_0 d s &\leftrightarrow &([a_1, a_0], s).
\end{array}
\end{equation}
On the left hand side, $\tilde{a}= \tilde{a}_1+ \tilde{a}_0 d s$ is an
$A$-path of $T^*(M \times \R)$ over the base path $\gamma= (\gamma_0,
\gamma_1)$ in $M\times\R$,
while on the right hand side we have a pair consisting of an $A$-path $a=
(a_1, a_0)$ of $T^*M
\oplus_M \R$ over the base path $\gamma_1$ in $M$, together with a real
number $s$. The explicit
formulas are:
\[ a_i(t) =e^{-\gamma_0(t)} \tilde{a}_i (t), \;\;\;\;\;
(i=0,1),\quad  s=\gamma_0(0) \] for computing the right hand side
from the left one, while for the other direction:
\[ \begin{split}
\tilde{a}_i(t) &=e^{\gamma_0(t)} a_i(t),\\
\gamma_0(t) &= -\int_0^t i(E)a_1(t)d t + s.
\end{split} \]
\end{remark}

\section{The Poisson case I: general discussion}
\label{Poisson-I}

\begin{num}\underline{{\bf The main theorem of the section.}}\rm \
In this section (as well as in the next two) we look at what the Jacobi
point of view brings new into Poisson geometry. In other words, we start
with a Poisson manifold $(P, \Lambda)$, and we view it both as a Poisson manifold, as well as
a Jacobi manifold with trivial Reeb vector field. Then $(P, \Lambda)$
will have two associated groupoids: the symplectic monodromy
groupoid $\Sigma_s(P)$ (with the role of integrating the Poisson structure),
and the contact monodromy groupoid $\Sigma_c(P)$ (with the role of integrating the Jacobi structure).
The aim of this section is to describe the relation
between the two.
Emphasize here that the smoothness of the one does not imply the
smoothness of the other, i.e. $P$ can be integrable as Poisson
manifold without being integrable as
Jacobi, or the other way around (see the examples).

Given $(P, \Lambda)$, we will define a bundle of groups
$\mathcal{P}_{\Lambda}$ over $P$,
where each fiber $\mathcal{P}_{\Lambda, x}$ is an additive subgroup of $\R$.
When $P$ is integrable
as a Poisson manifold, i.e. when $\Sigma_{s}(P)$ is symplectic Lie groupoid,
with symplectic form
denoted by $\Omega$, then $\mathcal{P}_{\Lambda, x}$ can be described as the
group of periods
of $\Omega|_{\alpha^{-1}(x)}$:
\begin{equation}
\label{periods}
Per(\Omega|_{\alpha^{-1}(x)})= \{ \int_{g} \Omega: [g] \in
\pi_2(\alpha^{-1}(x))\}.
\end{equation}
It is remarkable that these groups can be defined without any integrability
assumption. We also consider the quotient $\mathcal{G}_{\Lambda}$ of $P\times\mathbb{R}$ by $\mathcal{P}_{\Lambda}$,
i.e. the bundle of groups over $P$ whose fiber above $x\in P$ is $\R/\mathcal{P}_{\Lambda, x}$.
The main result to be discussed in this section is the following improvement of Main Theorem \ref{main2}
stated in the introduction.
\end{num}

\begin{theorem}
\label{poisson1}
For any Poisson manifold $(P, \Lambda)$, there is a short exact sequence of
topological groupoids
\[ 1\rmap \mathcal{G}_{\Lambda}\rmap \Sigma_{c}(P)\stackrel{\pi}{\rmap}
\Sigma_{s}(P)\rmap 1 .\] If $P$ is integrable as a Poisson
manifold, then the following are equivalent:
\begin{enumerate}
\item[(i)] $P$ is integrable as a Jacobi manifold.
\item[(ii)] $\mathcal{P}_{\Lambda}$ is locally uniformly discrete.
\item[(iii)] $\mathcal{P}_{\Lambda}$ is smooth.
\item[(iv)] $\mathcal{G}_{\Lambda}$ is smooth.
\end{enumerate}
Moreover, in this case the symplectic form $\Omega$ of $\Sigma_s(P)$ and the
contact form
$\theta$ of $\Sigma_c(P)$ are related by
\[ \pi^*\Omega= d\theta ,\]
and the Reeb vector field of $\Sigma_{c}(P)$ is
\[ R= \frac{d}{ds}, \]
the infinitesimal generator of the action of $\mathcal{G}_{\Lambda}$ on
$\Sigma_{c}(P)$
(or, equivalently, of the induced action of $\R$ via the projection $\R\rmap
\mathcal{G}_{\Lambda}$).
\end{theorem}

We mention here that, as
for the $G(A)$'s, when talking about
the smoothness of $\mathcal{P}_{\Lambda}$ or $\mathcal{G}_{\Lambda}$ we
refer to the natural smooth structures,
i.e. $\mathcal{P}_{\Lambda}$ as a subspace of $P\times\R$, and
$\mathcal{G}_{\Lambda}$ as a quotient of it. In particular,
there is at most one such smooth structure.
Also, the condition that $\mathcal{P}_{\Lambda}$ is locally uniformly discrete means that,
for $x\in P$, the distance between the zero element and
$\mathcal{P}_{\Lambda, y}-\{0_y\}$ is bounded from below by a positive number, for $y$ in
a neighborhood of $x$. \\

\begin{num}\label{describe-as-bdry}\underline{{\bf Monodromy maps.}}\rm \
To define the groups $\mathcal{P}_{\Lambda}$ in the nonintegrable
case (and to proceed with the proof of the theorem), we need to
recall the construction of the monodromy map of an algebroid (at a
first reading, readers can restrict themselves to the integrable
case, and skip the general definition of $\mathcal{P}_{\Lambda}$).

So, let $A$ be a Lie algebroid over $P$. For $x\in P$ we denote by
$\mathfrak{g}_x(A)$ the kernel of the anchor at $x$,
and we call it the isotropy Lie algebra of $A$ at $x$. The Lie bracket can be restricted to
this kernel, and
this shows that $\mathfrak{g}_x(A)$ is indeed a Lie algebra. As for any Lie
algebroid, we can form the
associated groupoid $G(\mathfrak{g}_x(A))$, which is a group since the base
is a point. This is nothing but
the unique simply connected Lie group integrating
$\mathfrak{g}_x(A)$, viewed as $\mathfrak{g}_x(A)$-homotopy classes
of paths $a_1: I\rmap \mathfrak{g}_x(A)$. Also, the image of $\rho$ defines
a singular foliation on $P$, whose leaves are the orbits of $A$.

\begin{defi}
Let $A$ be a Lie algebroid over $P$, let $x\in P$, and we denote by
$L_x$ the orbit through $x$. The monodromy map at $x$,
\[ \partial_{A}: \pi_2(L, x)\rmap G(\mathfrak{g}_x(A)) ,\]
associates to the homotopy class of $\gamma :I\times I\rmap L$
(representing an element in $\pi_2(L, x)$) the
class of a (any) $\mathfrak{g}_x$-path $a_1:I\rmap \mathfrak{g}_x$
with the property that there exists an
$A$-homotopy $a(\epsilon, t)$ (cf. Definition \ref{homotopy-def}), sitting over
$\gamma(\epsilon, t)$, and
which connects the zero path (i.e. $a(0, t)= 0$) to $a_1$.
\end{defi}

The image of $\partial_{A}$ is called the (extended) monodromy group of
$A$ at $x$, and these groups are central for understanding the integrability of $A$.
In particular, $A$ is integrable if and only if these groups are discrete,
locally uniformly with respect to $x$. For this, and more details (e.g.
to see that $a$ above can be chosen,
and that $[a_1]$ does not depend on the choice of $a$, but only on the homotopy
class of $\gamma$), we refer to \cite{marius}. We describe here what happens in
the integrable case, when these
constructions become more transparent. Then $G(A)$ is smooth,
and the isotropy group $G_x(A)$ (i.e. arrows of $G(A)$ that start
and end at $x$)
will be a Lie group that integrates the Lie algebra $\mathfrak{g}_x(A)$.
It may be different from $G(\mathfrak{g}_x(A))$, and the main reason is that
it may fail to be simply connected. Applying the homotopy long exact sequence
associated the fibration $\beta: \alpha^{-1}(x)\rmap L$, we get a surjective
boundary map
\[ \partial: \pi_2(L)\rmap \pi_1(G_x(A)) ,\]
and this is basically $\partial_A$ after one views $\pi_1(G_x(A))$ as a
subgroup of $G(\mathfrak{g}_x(A))$ (which is naturally the case since $G_x(A)$
integrates $\mathfrak{g}_x(A)$).
\end{num}

\begin{num}\underline{{\bf The period groups in the general case.}}\rm \
\label{general-m}
When $(P,\Lambda)$ is a Poisson manifold, we denote by
\[ \partial_s : \pi_2(L, x)\rmap G(\mathfrak{g}_{x}(T^*P)) ,\]
the monodromy map associated to $T^*P$, and, similarly, by $\partial_{c}$
the one associated to $T^*P\oplus \R$ (the algebroid associated to $P$
viewed as a Jacobi manifold). Note that $L$ is the symplectic leaf through $x$,
whose tangent space is spanned by the Hamiltonian vectors $X_f=
\Lambda^{\sharp}(df)= \{f, \cdot\}$,
and with the symplectic form
\[ \omega_{L}(X_f, X_g)= \{f, g\} .\]
\end{num}

\begin{defi}
Given a Poisson manifold $(P,\Lambda)$, $x\in P$, we define the period group
of $\Lambda$ at $x$,
\[ \mathcal{P}_{\Lambda, x}= \{ \int_{\gamma} \omega_L : [\gamma] \in
\pi_2(L,x) \; \text{and} \; \partial_s\gamma=1_x\}\subset \R,\]
and we define the period bundle $\mathcal{P}_{\Lambda}$ of
$\Lambda$ whose
fiber at $x$ is $\mathcal{P}_{\Lambda, x}$, and the structural group bundle
$\mathcal{G}_{\Lambda}$ of $\Lambda$ whose
fiber at $x$ is $\R/\mathcal{P}_{\Lambda, x}$.
\end{defi}

As promised, in the integrable case we have

\begin{lemma}
\label{expected}
If $P$ is integrable as a Poisson manifold, then $\mathcal{P}_{\Lambda, x}$
coincides with the group of
periods (\ref{periods}) of the restriction of the symplectic form of
$\Sigma_s(P)$ to
the $\alpha$-fiber at $x$.
\end{lemma}

\begin{proof}
Since $T^*P$ is integrable, we can use the description of $\partial_s$ as
the boundary of the homotopy long exact sequence (see the end of subsection \ref{describe-as-bdry}).
 We deduce that
\[ \mathcal{P}_{\Lambda, x}= \{ \int_{\beta_{*}(u)}\omega: u\in
\pi_2(\alpha^{-1}(x))\}.\]
On the other hand, since $\beta^{*}\omega= -\Omega$ on $\alpha^{-1}(x)$
($\beta$ is
anti-Poisson), we have
\[ \int_{\beta_{*}(u)}\omega= \int_{u}\beta^{*}\omega= -\int_{u}\Omega ,\]
and the lemma follows.
\end{proof}

Next, the two monodromy maps $\partial_{s}$ and $\partial_{c}$ are related as follows:

\begin{lemma}
\label{partial}
Let $P$ be a Poisson manifold and $x\in P$.
Denote by $\mathfrak{g}^{s}_{x}$ the isotropy Lie algebra at $x$ of $T^*P$
and by $\mathfrak{g}^{c}_{x}$ the isotropy Lie algebra at $x$ of $T^*P\oplus \R$.
Then
\[ G(\mathfrak{g}^{c}_{x})\cong G(\mathfrak{g}^{s}_{x})\times \R ,\]
and, after this identification, the monodromy maps $\partial_s$ and
$\partial_c$
of $T^*P$, and $T^*P\oplus \R$, respectively, are related by:
\[ \partial_c \gamma = ( \partial_s \gamma, - \int_{\gamma}
\omega_L), \]
for every $[\gamma] \in \pi_2(L, x)$.
\end{lemma}

\begin{proof}
The first part follows from the remark that $\mathfrak{g}^{c}_{x}=
\mathfrak{g}^{s}_{x}\times \R$, which is clear at the level of
vector spaces, and, at the level of Lie algebras, it follows
immediately from the formulas defining the Lie brackets of $T^*P$
and $T^*P\oplus\R$. Consider now $[\gamma] \in \pi_2(L)$, let
$\tilde{a}(\epsilon, t)$ be an $A$-homotopy over $\gamma$
connecting the trivial path with $\tilde{a}_1= (a_1, u_1)$ (so
that $\partial_c([\gamma])= [\tilde{a}_1]$), and let $b$ be the
solution of the equation (\ref{homotopy}). Write $a= (a, u)$ and,
similarly, $b= (b, v)$. The first component of the equation
\eqref{homotopy} tells us that $a$ is a homotopy between the zero
path and $a_1$, hence $\partial_c(\gamma)= (\partial_{s}(\gamma),
[u_1])$. On the other hand, $G(\R)\cong \R$, where the homotopy
class of an $\R$-path $u_1$ is identified (homotopic) to the
number $\int_{0}^{1} u_1$. We now look at the $\R$-component of
the equation \eqref{homotopy}, which gives us
\[\partial_t v -\partial_{\epsilon} u = \Lambda (a, b). \]
Since
\[ \Lambda^{\sharp} (a) = \frac{d\gamma }{d t} ,\;\;\;\;\;  \Lambda^{\sharp} (b)
= \frac{d\gamma}{d \epsilon} , \]
and $\gamma$ stays entirely in the leaf $L$, we have
\[ \partial_t v - \partial_{\epsilon} u = \omega_L ( \frac{d\gamma }{d t} ,
\frac{d\gamma}{d \epsilon}). \]
So
\[
\begin{split}
\int_0^1 \int_0^1 \omega_L( \frac{d\gamma }{d t} ,
\frac{d\gamma}{d \epsilon} )d \epsilon dt &= \int_0^1 \int_0^1
\partial_{t} v d \epsilon  d t -\int_0^1 \int_0^1 \partial_{\epsilon} u d \epsilon d
t \\
&= \int_0^1 (v(\epsilon , 1)- v(\epsilon , 0))d \epsilon- \int_0^1
( u(1,t)-u(0,t)) d t
\\
&=-\int_0^1 u(1, t) d t,
\end{split}
\]
i.e.  $\int_{\gamma} \omega_L = -\int_0^1 u(1,t) d t$. This proves
the lemma.
\end{proof}

\begin{num}\underline{{\bf Proof of Theorem \ref{poisson1}.}}\rm \
We now proceed with the proof of the theorem. The projection
$T^*P\oplus \mathbb{R}\rmap T^*P$ is an algebroid morphism, hence it induces
a groupoid morphism
\[ \pi: \Sigma_{c}(P)\rmap \Sigma_s(P) \]
which sends an $A$-path $(a, u)$ of $T^*P\oplus \mathbb{R}$ into the
$A$-path $a$ of $T^*P$. This is clearly surjective, and we denote by $\mathcal{G}$
its kernel. We will show that $\mathcal{G}= \mathcal{G}_{\Lambda}$.
Recall \cite{marius} that, for any algebroid $A$ over $P$, the isotropy group
at $x$ of the monodromy groupoid of $A$ (denoted by $G_x(A)$) has $\pi_0(G_x(A))$ isomorphic to
$\pi_{1}(L)$, and the connected component of the identity in $G_x(A)$ is
\[ (G_x(A))^0= G(\mathfrak{g}_x(A))/ Im(\partial_{A}) ,\]
independently of the integrability of $A$. Applying this to our algebroids, we see that
\[ \mathcal{G}_{x}= Ker(\pi: G(\mathfrak{g}_{x}^{c})/Im(\partial_{c})\rmap G(\mathfrak{g}_{x}^{s})/Im(\partial_{s})),\]
and, using the previous lemma, this is precisely $\mathbb{R}/\mathcal{P}_{\Lambda, x}$.
This proves the exact sequence in the theorem. Note that the inclusion of $\mathcal{G}_{\Lambda}$ into $\Sigma_{c}(P)$
sends the class of the real number $\lambda$ into the $A$-homotopy class of the path $(0, \lambda)$ of $T^*P\oplus \mathbb{R}$.

We now prove the equivalence of (i)-(iv) in the Theorem.
We first show that (i) is equivalent to (ii). By the general result of \cite{marius}, (i) is equivalent
to the groups $Im(\partial_{c, x})$ being locally uniformly discrete. This means that, if $(x_i)$ is a sequence in $P$
converging to $x$, and $[\gamma_i]\in\pi_2(L,x_i)$ satisfies
\begin{equation}
\label{lim}
\lim_{n\to + \infty} \text{distance} ((\partial_c (\gamma_i), 0))=0,
\end{equation}
then $\partial_c (\gamma_i)= 0$ for $i$ large enough. On the other hand, the similar condition
for $\partial_{s}$ is satisfied since $T^*P$ is integrable. Hence, using
Lemma \ref{partial}, the condition becomes: if $(x_i)$ is a sequence in $P$
converging to $x$, and $[\gamma_i]\in\pi_2(L,x_i)$ satisfies
\[ \int_{\gamma} \omega_{L_{x_i}}\rmap 0,\]
then these integrals must vanish for $i$ large enough. I.e.,
$\mathcal{P}_{\Lambda}$ must be locally uniformly discrete. Next,
(i) implies (iv) because $\pi$ will be a submersion and
$\mathcal{G}_{\Lambda}= \pi^{-1}(P)$ ($P$ is embedded in
$\Sigma_s(P)$ as the space of identity elements). Similarly, (iv)
implies (iii) because $\mathcal{P}_{\Lambda}$ is $\pi_{0}^{-1}(P)$
where the projection $\pi_{0}: P\times\mathbb{R}\rmap
\mathcal{G}_{\Lambda}$ is a submersion. Assume now (iii), and
prove (ii). The condition (iii) implies that
$\mathcal{P}_{\Lambda_x}$ is a smooth submanifold of $\mathbb{R}$.
But $\mathcal{P}_{\Lambda_x}$ is at most countable, since it is a
quotient of the second homotopy group $\pi_2(L)$ of the leaf
through $x$. It follows that $\mathcal{P}_{\Lambda_x}$ is discrete
and the projection from $\mathcal{P}_{\Lambda}$ into $P$ is a
local diffeomorphism. This implies (ii).

We are left with proving the last part of the theorem.
Using the
correspondence we have established in \eqref{corr} and the formula
of the symplectic form in the path space \cite{cafe}, we have
\begin{equation} \label{basic}
\mathcal{L}_{\frac{d}{d s}}
\theta =0, \;\;\;\;\; i(\frac{d}{d s}) \theta =1,
\end{equation}
i.e. $\frac{d}{d s}$ is the Reeb vector field of
$\theta$. From these formulas it also follows that $d\theta$ is a
basic form, i.e. $d \theta =\pi^*\omega$ for some 2-form $\omega$ on $\Gamma_s(P)$.
Since $\theta$ is a contact form and it is multiplicative, it follows that
$\omega$ is a symplectic form on $\Sigma_{s}(P)$ which is multiplicative.
Since the source map of $\Sigma_{c}(P)$ is Jacobi, it follows that
the source map $\alpha: (\Sigma_{s}(P), \omega)\rmap P$ is Poisson.
By uniqueness of the symplectic groupoid integrating the Poisson manifold $P$,
we must have $\omega=\Omega$.
\end{num}

\section{The Poisson case II: Application to prequantization}

\begin{num}\underline{{\bf The main theorem of the section.}}\rm \
Recall that a prequantization of a symplectic manifold $(S,
\omega)$ is a complex line bundle $L$ together with a connection
$\nabla$ so that $\omega$ represents the first Chern class $c_1(L,
\nabla)$. Equivalently, this is the same as a principal
$S^1$-bundle $\pi: \tilde{S}\rmap S$ together with a connection
1-form $\theta\in\Omega^1(\tilde{S})$ with the property that
$\pi^*\omega=d\theta$. Kostant's theorem (sometimes also
attributed to Kobayashi \cite{blair} or to Souriau) says that this
is possible if and only if $\omega$ is integral. Similarly,
Weinstein and Xu have introduced the notion of prequantization of
symplectic groupoids, with the aim of quantizing Poisson manifolds
``all at once" \cite{prequan}. More precisely:
\end{num}
\begin{defi} One calls prequantization of the symplectic groupoid
$(\Sigma, \Omega)$
any Lie groupoid extension of $\Sigma$ by the trivial bundle of Lie groups
$S^1$,
\[ 1\rmap S^1\rmap \tilde{\Sigma} \stackrel{\pi}{\rmap} \Sigma \rmap 1,\]
(and this makes $\tilde{\Sigma}$ into a principal $S^1$-bundle over
$\Sigma$),
together with a connection 1-form $\theta\in\Omega^1(\tilde{\Sigma})$
which is multiplicative and
satisfies $\pi^*\Omega=d \theta$.
\end{defi}

When saying ``the trivial bundle of $S^1$'s'', we really mean
triviality in the sense of extensions, i.e., besides the
triviality as a bundle, the action of $\Sigma$ on the bundle must
be trivial too. That simply means that $gz= zg$ for all $g\in
\tilde{\Sigma}$ and $z\in S^1$. In particular, there is no
ambiguity when talking about $\tilde{\Sigma}$ as a principal
$S^1$-bundle over $\Sigma$. This corresponds to the ``no-holonomy
above identity elements'' condition that appears in
\cite{prequan} (where uniqueness is proven).

In this section we show that, for a Poisson manifold $P$,
$\Sigma_{c}(P)$ is intimately related to prequantizing
the symplectic groupoid $\Sigma_{s}(P)$.
More precisely, we will prove the following result which is an
extension of Kostant's theorem to symplectic groupoids, and
an improvement of Main Theorem \ref{main3} stated in the introduction.

\begin{theorem}
\label{prqz}
Let $P$ be an integrable Poisson manifold, with associated
symplectic groupoid $\Sigma_s(P)$. The following are equivalent:
\begin{enumerate}
\item[(i)] $\Sigma_s(P)$ is prequantizable.
\item[(ii)]  $\mathcal{P}_{\Lambda}\subset P\times \mathbb{Z}$. 
\end{enumerate}
Moreover,  if $\Sigma_s(P)$ is Hausdroff, the
conditions above are also equivalent to 
\begin{enumerate}
\item[(iii)] $\Omega$ is integral.
\end{enumerate}
Finally, in the prequantizable case, the prequantization $\tilde{\Sigma}$ 
together with the
connection $1$-form becomes a contact groupoid which is a quotient of
$(\Sigma_{c}(P), \theta)$, and $\tilde{\Sigma}$ is Hausdorff if $\Sigma_{c}(P)$ is.
\end{theorem}

\begin{num}\underline{{\bf Proof of Theorem \ref{prqz}}}\rm \
The implications (i)$\Longrightarrow$ (ii) follows from the fact that the $s$-fibers of a prequantization groupoid are classical prequantizing bundles for the $s$-fibers of the symplectic groupoid. Also, in the Hausdorff case, the implications (i)$\Longrightarrow$ (iii)$\Longrightarrow$ (ii)
are clear.
Assume now (ii), and we prove (i) and the last part of the theorem. First of all, it follows that
$\mathcal{P}_{\Lambda}$ is uniformly discrete, hence, by Theorem
\ref{poisson1}, $\Sigma_c(P)$ is a Lie groupoid.

We put $\tilde{\Sigma}:= \Sigma_c(P)/\Z$, where the action of
$\mathbb{Z}$ is the one induced by the action of $\mathbb{R}$ (see
Theorem \ref{poisson1}). From the hypothesis and Theorem
\ref{poisson1} it also follows that $P\times
\mathbb{Z}/\mathcal{P}_{\Lambda}$ is a smooth (\'etale) sub-bundle
of $\mathcal{G}_{\lambda}$. Since $\Sigma_c(P)$ is a principal
$\mathcal{G}_{\Lambda}$-bundle over $\Sigma_{s}(P)$, it follows
that $\tilde{\Sigma}$ is smooth, and it is a principal
$\mathcal{G}_{\Lambda}/\mathbb{Z}= \mathbb{R}/\mathbb{Z}=
S^1$-bundle over $\Sigma_s(P)$. Denote by $\tilde{\pi}:
\tilde{\Sigma}\rmap \Sigma_s(P)$ the projection. This will be a
morphism of Lie groupoids, whose kernel is the trivial bundle of
groups with fiber $S^1$.

By \eqref{basic}, $\theta$ is $\R$-invariant, so it descends to a
1-form $\tilde{\theta}\in \Omega^1(\tilde{\Sigma})$ such that
$\tilde{\pi}^*\tilde{\theta}=\theta$. The same equations
\eqref{basic} imply the similar equations for $\tilde{\theta}$,
where the Reeb vector field will be the generator of the action of
$S^1$ on $\tilde{\Sigma}$. This shows that $\tilde{\theta}$ is a
connection 1-form on our principal bundle. Moreover, $\pi^*\Omega=
d\theta$ implies that $\tilde{\pi}^*\Omega= d\tilde{\theta}$.
Hence $(\tilde{\Sigma}, \tilde{\theta})$ is a prequantization of
our symplectic groupoid $\Sigma_s(P)$. By construction, it is a
quotient of $(\Sigma_{c}(P), \theta)$ and it inherits the contact
groupoid structure from $\Sigma_c(P)$. When $\Sigma_s(P)$ is
Hausdorff, so is $\tilde{\Sigma}$ as a $S^1$ principal bundle over
$\Sigma_s(P)$.  
\end{num}

\section{The Poisson case III: Relation to the Van Est map}

\begin{num}\underline{{\bf The main theorem of the section.}}\rm \
In this section we discuss the Poisson bivector from the point of view
of 2-cocycles. The main remark is that the algebroid
$T^*P\oplus\mathbb{R}$ (provided by the Jacobi point of view) is
made out of the algebroid $T^*P$ (provided by the Poisson point of
view) together with an extra-piece of data, namely a 2-cocycle.
This will show that most of our results belong to a more general
class of results, depending on closed 2-cocycles on algebroids. As
a consequence of our discussions, labelled here as the main
theorem of the section, we have the following result which gives
the precise conditions for when the relation between
$\Sigma_{c}(P)$ and $\Sigma_{s}(P)$ is the simplest possible. More
precisely, we will prove the following:
\end{num}

\begin{theorem}
\label{poisson2}
Let $(P, \Lambda)$ be a Poisson manifold, and assume that the symplectic monodromy groupoid
$\Sigma_{s}(P)$ is smooth and Hausdorff. Denote by $\Omega$ the symplectic form on $\Sigma_{s}(P)$.
Then the following are equivalent:
\begin{enumerate}
\item[(i)] $\Lambda$ is integrable as an algebroid cocycle.
\item[(ii)] $\Omega$ is exact.
\item[(iii)] $\mathcal{P}_{\Lambda}= 0$.
\end{enumerate}
Moreover, if $c$ integrates $\Lambda$, then $\Sigma_{c}(P)\cong \Sigma_{s}(P)\ltimes_{c}\R$.
\end{theorem}

This section is organized as a discussion around this theorem, which will provide
the precise definitions and the proof.

\begin{num}\underline{{\bf 2-cocycles and the Van Est map.}}\footnote{here we only give a brief outline
on algebroid cocycles, groupoid cocycles, and the Van Est map.
More details can be found in \cite{m-von-est}.}\rm \
\label{cocycles-sub} To see that 2-cocycles are at the heart of
$T^*P\oplus\mathbb{R}$ we remark that, if $[\cdot, \cdot]$ and
$\rho$ are the Lie bracket and the anchor of $T^*P$, then the
bracket of $T^*P\oplus\mathbb{R}$ is
\[
\begin{split}
[(\omega_1,\omega_0),(\eta_1,\eta_0)]  = & \big( [\omega_1,
\eta_1], \\ & \mathcal{L}_{ \rho
(\omega_1)}(\eta_0)-\mathcal{L}_{\rho
(\eta_1)}(\omega_0)+  \Lambda(\eta_1,\omega_1) \big)\\
\end{split}
\]
(and the anchor is $(\omega_1,\omega_0)\mapsto \rho(\omega_1)$).
More abstractly, given any algebroid $A$ and any section $\Lambda\in\Gamma(\Lambda^2A^*)$, the previous formula
defines a bracket on $A\oplus\mathbb{R}$. On the other hand, the spaces $\Gamma(\Lambda^pA^*)$ together define
an ``$A$-De Rham complex'', with the differential
\begin{eqnarray}\label{differential}
d_{A}(\Lambda)(X_1, \ldots , X_{p+1}) & = & \sum_{i<j}
(-1)^{i+j-1}\Lambda([X_i, X_j], X_1, \ldots , \hat{X_i}, \ldots ,
\hat{X_j}, \ldots X_{p+1})) \nonumber \\
 & + & \sum_{i=1}^{p+1}(-1)^{i}
\mathcal{L}_{X_i}(\Lambda(X_1, \ldots, \hat{X_i}, \ldots ,
X_{p+1})) 
\end{eqnarray}
(where $\Lambda\in \Gamma(\Lambda^pA^*)$ is arbitrary). 
The resulting complex is denoted by $\Omega^*(A)$, and the cohomology
is denoted by $H^*(A)$. For instance, when $A= TP$
one recovers the usual De Rham cohomology, and when $A= \mathfrak{g}$ is a Lie algebra
one recovers the Lie algebra cohomology.

For $\Lambda\in \Gamma(\Lambda^2A^*)$, the resulting structure on $A\oplus \mathbb{R}$ is a Lie algebroid structure
if and only if $\Lambda$ is a 2-cocycle, i.e. $d_{A}\Lambda= 0$. The resulting algebroid is denoted by $A\ltimes_{\Lambda}\mathbb{R}$.
It is not difficult to see that the isomorphism class of $A\ltimes_{\Lambda}\mathbb{R}$ only depends on the cohomology class of $\Lambda$.
Coming back to our Poisson manifold, we conclude that
\end{num}

\begin{lemma} Given a Poisson manifold $(P, \Lambda)$, the Poisson tensor is
a closed algebroid 2-cocycle for $T^*P$, and
\[ T^*P\oplus\mathbb{R}\cong T^*P\ltimes_{\Lambda}\mathbb{R} .\]
\end{lemma}

We now turn to the global picture, i.e. to cocycles on groupoids. Assume that $G$ is a Lie groupoid over $P$.
A differentiable $p$-cocycle on $G$ is a smooth function $c$ defined on the space of $p$-tuples
$(g_1, \ldots , g_p)$ of composable elements of $G$ (i.e. such that $g_1\ldots g_p$ is defined).
The differential of $c$ is the $(p+1)$-cocycle
\begin{eqnarray}
(dc)(g_{1}, \ldots , g_{p}, g_{p+1}) & = & c(g_{2}, \ldots , g_{p+1})  + \\
   &  & \sum_{i=1}^{p} (-1)^i c(g_{1}, \ldots , g_{i}g_{i+1}, \ldots , g_{p+1}) +  \\
   & & (-1)^{p+1}c(g_{1}, . . . , g_{p}) .
\end{eqnarray}
We will only work with normalized cocycles, i.e. with the property that $c(g_1, \ldots , g_p)$ $=0$
whenever one of the entries is an identity arrow. The resulting
complex is denoted by $C^{*}_{diff}(G)$, and
the cohomology is denoted by $H^{*}_{diff}(G)$.

As in the algebroid case, any (normalized) 2-cocycle $c$ induces a Lie groupoid structure on $G\times\mathbb{R}$,
with $(g_1, \lambda_1)(g_2, \lambda_2)= (g_1g_2, \lambda_1+ \lambda_2+ c(g_1, g_2))$.
We denote this groupoid by $G\ltimes_{c}\mathbb{R}$. Again, the isomorphism class of $G\ltimes_{c}\mathbb{R}$ only
depends on the cohomology class of $c$.

Finally, by differentiation, groupoid cocycles on $G$ induce algebroid cocycles on the
Lie algebroid $A$ of $G$ (see \cite{m-von-est}). This construction induces a
map of complexes $VE: C^{*}_{diff}(G)\rmap \Omega^*(A)$, hence also a map in cohomology
\[ VE: H^{*}_{diff}(G)\rmap H^{*}(A) ,\]
known as the Van Est map. Elements in the image of the Van Est map are called integrable
algebroid cocycles (note that this property only depends on the cohomology class of a cocycle).

\begin{defi}
Given a Poisson manifold $(P, \Lambda)$, we say that $\Lambda$ is integrable
if the Poisson manifold $P$ is integrable, and $\Lambda$ comes from a differentiable groupoid 2-cocycle $c$
on $\Sigma_s(P)$. In this case we also say that $c$ integrates $\Lambda$.
\end{defi}

In general, the algebroid associated to an integrable
algebroid 2-cocycle is integrable. More precisely, if
$c$ is a groupoid 2-cocycle on a Lie groupoid $G$ and $\Lambda$ is the induced algebroid 2-cocycle,
it follows that the Lie algebroid of $G\ltimes_{c}\mathbb{R}$ is precisely $A\ltimes_{\Lambda}\mathbb{R}$
\cite{m-von-est}.

\begin{num}\underline{{\bf Proof of Theorem \ref{poisson2}.}}\rm \
Clearly (ii) implies (iii), and we now prove the converse.
We use Theorem \ref{poisson1} and the fact that
$\G_x=\R$ (the hypothesis). It follows that
 $\pi: \Sigma_c(P) \rmap \Gamma_s(P)$ is a principal
$\R$-bundle, with connection form $\theta$ and curvature $\Omega$.
Since $\R$ is contractible, $\pi^*$ induces an isomorphism from
$H^2(\Gamma_s(P))$ to $H^2(\Sigma_c(P))$. So $[\pi^* \Omega]=[d
\theta]$ implies that $[\Omega]=0$, i.e. $\Omega$ is exact.

To prove that (i) is equivalent to (iii), we use the following
characterization of integrable algebroid cocycles from \cite{m-von-est}.
Assume that $\nu$ is an algebroid 2-cocycle on $A$, and that $G$ is an
$\alpha$-simply connected Hausdorff Lie groupoid integrating $A$.
Since each fiber $A_x$ is identified with the tangent space
of $\alpha^{-1}(x)$ at the identity element at $x$, using right translations,
$\nu$ will induce a 2-form on $\alpha^{-1}(x)$, call it $\omega_{\nu, x}$.
Then, $\nu$ is in the image of the Van Est map if and only if $Per(\omega_{\nu, x})= 0$
for all $x\in P$. In our case, note that $\omega_{\Lambda, x}= \Omega|_{\alpha^{-1}(x)}$,
hence the equivalence of (i) and (iii) follows.
\end{num}

\begin{corollary}
For a compact s-simply connected contact groupoid of an integrable
Poisson manifold, the Reeb vector field always has at least one
closed orbit.
\end{corollary}
\begin{proof}  If the Reeb vector field has no closed orbit, then according to  Theorem \ref{poisson1}, $\G_\Lambda$ is the
trivial $\R$ bundle. From Theorem \ref{poisson2}, one has
$\Sigma_c(P)\cong \Sigma_s(P)\times \R$, which contradicts with
the compactness. \end{proof}

\begin{corollary}
If every symplectic leaf in an integrable Poisson manifold $P$ has
exact symplectic form, then the symplectic form $\Omega$ of
$\Sigma_s(P)$ is also exact.
\end{corollary}
\begin{proof} It is a direct conclusion from Theorem
\ref{poisson2} and Lemma \ref{expected}.
\end{proof}

\section{Special cases and examples}

\begin{ex}[\underline{Symplectic manifolds}]
\label{example1}\rm\
Let $(S, \omega)$ be a symplectic manifold, and assume for simplicity that
$S$ is simply connected (in general, one has to replace the pair groupoid below with the
homotopy groupoid). Then the symplectic groupoid of $S$ (viewed as a Poisson manifold) is
\[ \Gamma_s(S)=(S \times S, (\omega, -\omega)) ,\]
the pair groupoid (source and target are the  projections, and the multiplication is
$(x, y)(y, z)= (x, z)$), endowed with the symplectic form $pr_{1}^{*}\omega- pr_{2}^{*}\omega$.
In this case $\mathcal{P}_{\Lambda}$ is a trivial group bundle
over $S$ with fiber $Per(\omega)$, hence, by Theorem \ref{poisson1}, $S$ is integrable as a Jacobi
manifold  if and only if
$Per(\omega)$ is a discrete group, i.e. $Per(\omega)= a\mathbb{Z}$ for some real number $a$.
The interesting case is when $a\neq 0$. In this case $\omega_{a}= \frac{1}{a}\omega$
is integral, hence we find a principal $S^1$-bundle $\pi: \tilde{S}\rmap S$ and a connection 1-form $\theta_{a}\in\Omega^1(\tilde{S})$
so that $\pi^*\omega_a= d\theta_a$. The gauge groupoid of $\tilde{S}$ is $\tilde{S}\otimes_{S^1}\tilde{S}$, the quotient
of the pair groupoid $\tilde{S}\times \tilde{S}$ by the diagonal action of $S^1$ (hence the base manifold is $\tilde{S}/S^1= S$).
Moreover, the 1-form $(\theta_a, -\theta_a)$ is basic, hence induces a 1-form $\bar{\theta}_a$ on $\tilde{S}\otimes_{S^1}\tilde{S}$,
and this makes the gauge groupoid into a contact groupoid. It is not difficult to see that
\[ \Sigma_c(S)= (\tilde{S}\otimes_{S^1}\tilde{S}, a\bar{\theta}_a, 1).\]

\end{ex}

\begin{ex}[\underline{Contact manifolds}]\rm\
Let $(M, \theta)$ be a contact manifold, and assume for simplicity that
$M$ is simply connected (as in the previous example, in the general case one has to replace the pair groupoid by the homotopy groupoid).
Then the contact groupoid of $M$ (viewed as a Jacobi manifold) can be described as follows. Consider
the product $M\times M\times \mathbb{R}$ of the pair groupoid
with $\mathbb{R}$. The 1-form
\[ \theta= -(\exp \circ p_3) p_1^* \theta_0+p_2^* \theta_0 \]
where $p_i$, $1\leq i \leq 3$, is the projection on the $i$-th component, will be a contact form,
and it will be multiplicative with respect to the Reeb cocycle $r= \exp\circ p_3$. With these,
\[ \Sigma_{c}(M)= (M\times M\times \mathbb{R}, \theta, r) .\]
\end{ex}

\begin{ex}[\underline{Vector fields}] \rm\
Given a vector field $X$ on a manifold $M$, one can view $(M, X)$ as a Jacobi manifold with
vanishing bivector.
Note that the orbits of the associated Lie algebroid $T^*M\oplus \mathbb{R}$
are precisely the orbits of $X$ hence, since they are $1$-dimensional, it follows that
$M$ is integrable as a Jacobi manifold. Hence the associated contact groupoid,
denoted here by $\Sigma_{c}(M, X)$, is smooth. Let us describe $\Sigma_{c}(M, X)$ in more detail.
First, let us mention two other simpler groupoids which are associated to $X$.

1. The flow of $X$, $\mathcal{D}(X)$, is probably the best known example of Lie groupoid.
One has $\mathcal{D}(X)\subset M\times \R$ as the domain of the local flow $\phi^{t}_{X}$ of $X$,
consisting of pairs $(x, t)$ with the property that $\phi^{t}_{X}(x)$ is defined. The elements
$(x, t)\in \mathcal{D}(X)$ are viewed as arrows from $\phi^{t}_{X}(x)$ into $x$, and the composition
is given by the rule $\phi_{X}^{t}\phi_{X}^{s}= \phi_{X}^{t+s}$.

2. In general, for any finite dimensional vector space $V$ and any vector $v\in V$,
one has an associated Lie algebra $\mathfrak{g}(v)$, which is $V^*$ endowed with the bracket
\[ [\alpha, \beta]= -\alpha(v)\beta+ \beta(v)\alpha .\]
The associated simply connected Lie group, denoted by $\mathcal{G}(v)$, can be described as follows:
\[ \mathcal{G}(v)= \{ \lambda\in V^*: \phi_{\lambda}:= Id_{V}+ \lambda v\in Aut^{+}(V)\}\]
(where $Aut^{+}(V)$ is the group of orientation preserving automorphisms of $V$).
The product $\lambda\eta$ is defined by
\[ u\mapsto \lambda(u)+ \eta(u)+ \lambda(v)\eta(u) ,\]
i.e. so that $\phi_{\lambda\eta}= \phi_{\lambda}\phi_{\eta}$.
Applying this to each $X_x\in T_xM$, we obtain a bundle of Lie groups
over $M$, denoted by $\mathcal{G}(X)$.

Note that $\mathcal{D}(X)$ acts on $\mathcal{G}(X)$: for each $(x, t)\in \mathcal{D}(X)$
viewed as an arrow from $y= \phi_{X}^{t}(x)$ into $x$,
$(d\phi_{X}^{t})_x: T_xM\rmap T_{y}M$ preserves $X$, hence it induces
a Lie group map from $\mathcal{G}(X_y)$ into $\mathcal{G}(X_x)$,
denoted by $\phi_{x, t}$. One then forms
the semi-direct product $\mathcal{G}(X)\rtimes \mathcal{D}(X)$, which consists of triples
\[ (\lambda, x, t)\ \ \text{with} \ \ (x, t)\in \mathcal{D}(X), \lambda\in \mathcal{G}(X_x),\]
such a triple is viewed as an arrow from $\phi^{t}_{X}(x)$ into $x$, and the composition is given by
\[ (\lambda, x, t)(\lambda', x', t')= (\lambda+ \phi_{x, t}(\lambda'), x, t+ t') .\]
With these, the contact groupoid is
\begin{equation}
\label{solution}
\Sigma_c(M, X)= \mathcal{G}(X)\rtimes \mathcal{D}(X) .
\end{equation}
In particular, this tells us that the period group of $X$ at $x$,
\[ Per_{x}(X)= \{ t: \phi_{X}^{t}(x)= x \} ,\]
acts on the Lie group $\mathcal{G}(X_x)$, and the isotropy group of $\Sigma_c(M, X)$
(i.e. arrows that start and end at $x$) is
\[ \Sigma_c(M, X)_x= \mathcal{G}(X_x)\rtimes Per_{x}(X).\]
To see (\ref{solution}), one looks at what happens at the infinitesimal level. The Lie
algebroid of $\mathcal{D}(X)$, denoted by $\mathbb{L}_{X}$,  is the trivial line bundle together with the bracket
$[f, g]_X:= -X(f)g+ X(g)f$ for $f, g\in \Gamma(\mathbb{L}_{X})= C^{\infty}(M)$, and the anchor is
given by multiplication by $X$. Moreover, the Lie algebras
$\mathfrak{g}(X_x)$ fit into a bundle of Lie algebras, denoted
 by
$\mathfrak{g}(X)$; this is $T^*M$, together with the Lie bracket on
1-forms $[\omega, \theta]_{X}:= -i_{X}(\omega\wedge\theta)$. There is an obvious exact sequence of
algebroids
\[ 0\rmap \mathfrak{g}(X) \rmap T^*M\oplus_M\R \rmap \mathbb{L}_{X}\rmap 0,\]
and, from the explicit formulas for the bracket on $T^*M\oplus \R$, it is not difficult to see
that $T^*M\oplus_M\R$ is a semi-direct product of Lie algebroids, where the action of $\mathbb{L}_{X}$
on $\Gamma(\mathfrak{g}(X))= \Omega^1(M)$ is the Lie derivative with respect to $X$. Passing to the global picture,
we find (\ref{solution}).
\end{ex}

\begin{ex}[\underline{Homogeneous Poisson manifolds}]\rm\
Let $(P, \tilde{\Lambda})$ be a homogeneous Poisson manifold (see subsection \ref{hom-Ps}).
Note that the homogeneity equation $\mathcal{L}_{Z}\tilde{\Lambda}= -\tilde{\Lambda}$ can be reformulated
in terms of the cohomology complex of the algebroid $A= T^*P$ (see subsection \ref{cocycles-sub})
as $\Lambda= -d_{T^*M}(Z)$. In particular, Theorem \ref{poisson2} tells us that $\mathcal{P}_{\Lambda}= 0$.
We deduce the following:

\begin{corollary} Any integrable homogenous Poisson manifold is also Jacobi integrable,
and $\Sigma_{c}(P)\cong \Sigma_s(P)\times \R$.
\end{corollary}
\end{ex}

\begin{ex}[\underline{Conformal versions}]\rm\
Given a Jacobi manifold $(M, \Lambda, E)$, and a smooth nowhere vanishing function
$\tau$ on $M$, one defines the conformal transformation of $(\Lambda, E)$
by $\tau$ as the new Jacobi structure given by
\[ \Lambda_{\tau} = \tau \Lambda, \;\;\;\;\; E_{\tau}= \tau E+ \Lambda^{\sharp} (d\tau) .\]
One says that two Jacobi structures are conformal equivalent if they
are related by such a transformation; such an equivalence class
of Jacobi structures is called {\it conformal Jacobi structure}.
Of course, when restricted to contact manifolds, this becomes the usual
notion of conformal equivalence of contact forms (see e.g. \cite{blair}): $\theta$ and $\theta'$ are equivalent
if $\theta'= \tau \theta$ for some non-vanishing function $\tau$. Equivalently, this corresponds to the
fact that the contact distribution of $\theta$, $\mathcal{H}_{\theta}= Ker(\theta)$, coincides with
the one of $\theta'$.


Similarly, given a contact groupoid $\Sigma$ over $M$ with contact form $\theta$ and Reeb
function $r$, and a smooth nowhere vanishing function
$\tau$ on $M$,
\[ \theta_{\tau}= \alpha^*(\tau) \theta, \quad r_{\tau}= r+ \lg(\frac{\alpha^*\tau}{\beta^*\tau})\]
define a new contact form and Reeb function so that $(\Sigma, \theta_{\tau}, r_{\tau})$
is a contact groupoid. Exactly as before, $\Sigma$, together with an equivalence class
of a pair $(\theta, r)$, will be called a {\it conformal contact groupoid}.


We then have:

\begin{corollary} There is a 1-1 correspondence between conformal contact groupoids over $M$ which are $\alpha$-simply connected, and integrable conformal Jacobi structures on $M$.
\end{corollary}

This follows immediately from our results and the following two simple remarks:
\begin{enumerate}
\item[(i)] The integrability of Jacobi structures is stable under conformal equivalences, hence one
can talk about the integrability of conformal Jacobi structure.
\item[(ii)] If $f: (N, \Lambda^{'}, E^{'})\rmap (M, \Lambda, E)$
is a Jacobi map, and $\tau$ is a non-vanishing function
on $M$, then $f$ is Jacobi also as a map from $(N, \Lambda^{'}_{\tau\circ f}, E^{'}_{\tau\circ f})$ to $(M, \Lambda_{\tau}, E_{\tau})$.
\end{enumerate}
\end{ex}

\begin{ex}[\underline{Locally conformal versions}]\rm\
\label{loclly-conf}
Apparently, the contact groupoid equation
\[ m^* \theta = pr_2^*(e^{-r}) \cdot pr_1^* \theta + pr^*_2 \theta\]
is not symmetric. The ``mirror symmetry'' of the previous equation is
\[ m^* \theta '= pr_1^* \theta' + pr_1^*(e^{-r'})\cdot pr^*_2 \theta' ,\]
and this is obtained by the transformation $\theta'= e^r\theta$
and $r'= -r$. A more symmetric version of the equation is obtained
by choosing $\theta_0= e^{-\frac{r}{2}}\theta$, $r_0=
\frac{r}{2}$, for which we have
\[ m^* \theta_0 = pr_2^*(e^{-r_0}) \cdot pr_1^* \theta_0 + pr_1^*(e^{r_0})pr^*_2 \theta_0.\]
Of course, all these describe essentially the same contact
groupoid, and this is the point of view adapted in \cite{dazord0,
dazord}. The relation between all these descriptions comes from
the fact that $Ker(\theta)= Ker(\theta')= Ker(\theta_0)$. This is
only one of the several motivations for considering ``locally
conformal versions'' of the structures we have been discussing.

To formulate our conclusions, we recall here that a locally conformal Jacobi structure on $M$  is described by
\begin{enumerate}
\item[] \it{(LCJ1)} an open cover $\{U_i\}$ of $M$ together with Jacobi structures
$(\Lambda_i, R_i)$ on each open $U_i$, so that, on the overlaps $U_i\cap U_j$,
the restrictions of the two Jacobi structures are conformal equivalent via $\tau_{i, j}$,
and so that the $\tau_{i, j}$'s satisfy the cocycle condition $\tau_{i, j}\tau_{j, k}= \tau_{i, k}$ on $U_i\cap U_j\cap U_k$.
\end{enumerate}

Of course, different covers and local Jacobi structures can lead to the same locally conformal Jacobi structure, i.e.
one has to consider a certain equivalence relation. This is completely similar to the description
of vector bundles in terms (of equivalence classes!) of transition functions, and, as there,
there is an alternative global description:

\begin{enumerate}
\item[] \it{(LCJ2)} A (isomorphism class of a) line bundle $L$ over $M$ together with a Lie algebra structure
$[\cdot, \cdot]$ on the space $\Gamma(L)$ of sections, so that $[\cdot, \cdot]$ is local.
(The isomorphisms are realized by bundle maps covering the identity
and inducing Lie algebra isomorphisms.)
\end{enumerate}

The global picture (LCJ2) is obtained by gluing the local data of
(LCJ1): one glues the trivial line bundles over $U_i$ using the
transition functions $\tau_{i, j}$, and then the local Lie
brackets defined on each $U_i$ by the Jacobi structures (see the
introduction) will fit together into a local Lie bracket on
$\Gamma(L)$.

Restricting to Jacobi structures coming from contact ones, we obtain what we will call here locally conformal contact structure,
and which are well known in the literature (often under various other names). Similar to the discussion above (and well known), such a structure
on $M$ is described by either
\begin{enumerate}
\item[] \it{(LCC1)} an open cover $\{U_i\}$ of $M$ together with contact forms $\theta_i$
on each $U_i$ and nowhere vanishing functions $\tau_{i, j}$ defined on the overlaps $U_i\cap U_j$,
such that $\theta_j= \tau_{i, j}\theta_i$.
\item[] \it{(LCC2)} a contact hyperplane $\mathcal{H}$, i.e. a smooth field of tangent hyperplanes $\mathcal{H}\subset TM$ so that, locally,
$\mathcal{H}$ is of type $Ker(\theta)$ for some (locally defined) contact 1-form $\theta$.
\end{enumerate}

Note that a conformal Jacobi structure is the same thing as a
locally conformal Jacobi structure with orientable line bundle,
and a conformal contact structure is the same thing as a locally
conformal contact structure whose contact hyperplane is induced by
a globally defined contact form. In particular, if $M$ is simply
connected, then ``locally conformal= conformal''.

Now, with our terminology, a locally conformal contact groupoid will be {\it a groupoid $\Sigma$ over $M$ together
with a contact hyperplane $\mathcal{H}$ with the property that $\mathcal{H}$ is compatible with the groupoid structure in the sense that
\begin{enumerate}
\item[(i)] The inversion $i: \Sigma \rightarrow \Sigma$ leaves
$\mathcal{H}$ invariant.
\item[(ii)] For all $X,Y\in \mathcal{H}$ for which $X \cdot Y= (dm)(X, Y)$ is defined,
we have: $X \cdot Y \in \mathcal{H}$.
\end{enumerate}
}

These have been introduced in \cite{dazord0, dazord} under the name of contact groupoids (see also \cite{ZaZh}). With this, our main theorem on
contact groupoids and the correspondence with Jacobi structures has a locally conformal version (and this completes the study
of \cite{dazord0, dazord}). There are various ways to see this. For instance,
Dazord shows that, if $(\Sigma, \mathcal{H})$ is a locally conformal contact groupoid so that $\mathcal{H}= Ker(d\theta)$ for some globally defined contact form
$\theta$, then $\theta$ can be choosen so that $(\Sigma, \theta, r)$
is a contact groupoid (for some uniquely defined multiplicative
function $r$), and two choices $\theta_1$ and $\theta_2$ define the
same locally conformal contact groupoid if and only if the associated
contact groupoids are conformal equivalent (see Proposition 4.1 in
\cite{dazord} or Appendix I in \cite{ZaZh}). In particular, if $M$ is simply connected and $\Sigma$ is $\alpha$-simply connected, then locally conformal structures on $\Sigma$ are the same thing as conformal structures on $\Sigma$. One can then pull-back $\Sigma$ to the
universal cover $\tilde{M}$, use our main result there, and show that it decends down to $M$ (this requires some care; in particular, the trivial line bundle over $\tilde{M}$ will descend to a possibly non-trivial line bundle over $M$). Alternatively, one can check that all our arguments, after suitable modifications, make sense in the locally conformal setting as well. For instance, if $M$ is locally conformal with associated bundle $L_{M}$, and if $(\Sigma, \mathcal{H})$ is a locally conformal groupoid, then:
\begin{enumerate}
\item[(i)] the Poissonization of $(M, L_{M})$ will be $L_{M}$ viewed as a manifold.
\item[(ii)] the Lie algebroid of $(M, L_{M})$ will be the jet bundle $J^1(L_M)$, with the bracket of two 1-jets given by the jet of the local Lie bracket on $\Gamma(L_M)$.
\item[(iii)] the symplectification of $(\Sigma, \mathcal{H})$ will be $L_{\Sigma}= T\Sigma/\mathcal{H}$, a symplectic groupoid over $L_{M}$.
\end{enumerate}
All these have been already explained in \cite{dazord}. A bit more care is needed when working with $A$-paths. Nevertheless, one can use a connection on $L_{M}$ to write the jet algebroid as $L_M\oplus \mathbb{R}$, so that the discussion from section \ref{section3} (where the Lie bracket, the $A$-paths, and the corresponding equations are all written componentwise) can be carried out in this setting. In particular, we have:

\begin{corollary} There is a 1-1 correspondence between locally conformal contact groupoids over $M$ which are $\alpha$-simply connected, and integrable locally conformal Jacobi structures on $M$.
\end{corollary}
\end{ex}

\begin{ex}[\underline{De-poissonization}]\rm\
Inverse to the ``poissonization process'' (Section
\ref{Poissonization and homogeneity}), one can obtain Jacobi-type
structures out of homogeneous Poisson manifolds. The result is
also a very good illustration of the different Jacobi-type
structures described in the previous two examples. More precisely,
given a homogeneous Poisson manifold $(P, \tilde{\Lambda})$ with
homogeneous vector field $Z$, assuming that $Z$ is nowhere zero
and that the set $\bar{P}$ of all trajectories of $Z$ admits a
manifold structure such that the projection $\pi: P\rmap \bar{P}$
is a submersion, then
\begin{enumerate}
\item[(i)] $\bar{P}$ has an induced locally conformal Jacobi structure.
\item[(ii)] if $Z$ is the infinitesimal generator of a free action of $\R_{+}$ on $P$, then
$\bar{P}$ has a canonically induced conformal Jacobi structure.
\item[(iii)] under the condition of (ii), if $\pi$ has a distinguished section, then $\bar{P}$ has
a distinguished Jacobi structure.
\end{enumerate}
Part (i) is explained in \cite{dazord2} (Corollary 2.5), while (ii) and (iii) are just simpler versions that
we explain here independently of (i).
For (iii), using the given section we identify $P$ with $\bar{P}\times \mathbb{R}_{+}$  with the action on
the second component. Then, completely similar to Section
\ref{Poissonization and homogeneity}, one gets an induced Jacobi
structure on $\bar{P}$. For (ii), since the fibers of $\pi$ are contractible, we can allways find globally defined sections,
and consider the induced Jacobi structures induced from (iii); different choices of sections produce conformal
equivalent Jacobi structures, i.e. it is only the conformal Jacobi structure that is independent of the choice of
the section.

Of course, starting with a symplectic manifold, the quotient will be a (locally conformal) contact manifold. When
looking at a homogeneous symplectic groupoid $\Sigma$ over $P$, one can show that the induced quotient $\bar{\Sigma}$ inherits
a groupoid structure, and it is a (locally conformal) contact groupoid over $\bar{P}$. Of course, the quotient of $\Sigma_{s}(P)$
will coincide with $\Sigma_c(\bar{P})$.
\end{ex}

\begin{ex}[\underline{Sphere bundles}]\rm\
A particular case of the previous example comes from duals of Lie algebroids. Relevant here is that, given a
vector bundle $A$ over $M$, there is a 1-1 correspondence between Poisson structures on $A^*$
and Lie algebroid structures on $A$ (see e.g. \cite{ana}). In particular, given a Lie algebroid $A$, the Poisson
manifold $A^*$ will be homogeneous with respect to the generator of the (fiberwise linear!)
action of $\mathbb{R}_{+}$: $t\cdot a= ta$. It follows that the sphere bundle
\[ S(A^*)= (A^*-\{0\})/ \mathbb{R}_{+} \]
has a conformal Jacobi structure. When $A$ comes equipped with a
metric, $S(A^*)$ takes the more familiar form consisting of
vectors of norm $1$, and we will have an induced Jacobi structure
on $S(A^*)$ (in (iii) of the previous example, one uses the obvious section induced by the
metric). This applies in particular to $A= TM$, when one obtains $S(T^*M)$ with
its canonical contact structure (or conformal contact if one does not fix a metric).

On the other hand, if $A$ comes from a Lie groupoid $G$
that we assume to be $\alpha$-simply connected, one knows that $T^*G$ is
naturally a groupoid over $A^*$ which, together with the canonical
symplectic form on the cotangent bundle, becomes a symplectic
groupoid (cf. e.g. \cite{dazord}). Of course, this is the symplectic monodromy groupoid of
$A^*$. Passing to sphere bundles,
$S(T^*G)$ together with the canonical (conformal) contact
structure and the induced groupoid structure, becomes the
conformal contact groupoid associated to $S(A^*)$. Interesting
particular cases are obtained when $A$ is a Lie algebra, or a
tangent bundle.
\end{ex}

\begin{ex}[\underline{2-dimensional Poisson case}]\rm\
It is known that every 2-dimensional Poisson manifold $P$ is
integrable \cite{marius2}. The main reason is that each
symplectic leaf of $P$ is either a point or 2 dimensional, and in the former case
the isotropy Lie algebra of $T^*P$ at $x$ is zero. These force that
all monodromy maps $\partial_{s}$ (see subsection \ref{general-m})
are zero, hence, by the main result of \cite{marius}, $T^*P$ must be integrable.
Note that the previous discussion also shows that $\mathcal{P}_{\Lambda, x}$ coincides
with the period group $Per(\omega_L)$, where $L$ is the symplectic leaf through $x$ and
$\omega_L$ is its symplectic form. Hence, Theorem \ref{poisson1} becomes

\begin{corollary}
A 2-dimensional Poisson manifold is integrable as a Jacobi
manifold if and only if $Per(\omega_L)$ are locally uniformly discrete.
\end{corollary}
\end{ex}

\begin{ex}[\underline{Non-integrable examples}] \rm\
Using the previous corollary or Example \ref{example1}, it is easy
to produce Poisson manifolds which are integrable as Poisson
manifolds but which are not integrable as Jacobi manifolds: it
suffices to consider a symplectic form whose period group is
dense. We now show that there are also Poisson manifolds which are
Jacobi integrable, but which are not Poisson integrable. This
exploits one of the examples of \cite{marius2} (see Example 3.8
there), which we now recall. Let $M_a=\R^3$, endowed with the
Poisson structure described by the bracket of the coordinate
functions $x^i$ as follows:
\[ \{x^2, x^3 \} = a x^1 ,\;\;\; \{x^3, x^1 \} = a x^2, \;\;\; \{x^1, x^2 \}= a x^3, \]
where $a=a(r)$ is a function depending only on the radius $r$, with the property
that $a(r)> 0$ for $r> 0$. The symplectic leaves of $M_a$ are the spheres $S_r$ centered at the
origin (including the degenerated sphere:  the origin itself), and the leafwise symplectic forms
are
\[ \omega_{r}= \frac{r^2}{a}(x^1dx^2dx^3 + x^2dx^3dx^1+ x^3dx^1dx^2) .\]
Central to the conclusion is the symplectic area function,
\[ A_a(r)= \int_{S_r} \omega_r= \frac{4\pi r}{a(r)} .\]
With these, \cite{marius2} shows that $M_a$ is Poisson integrable
if and only if either $A_{a}(r)$ is constant, or $A_{a}$ has no critical points and
$\lim_{r\to 0}A'_a(r) \neq 0$. One can do exactly the same type
of computations as in \cite{marius2}, but this time for the algebroid
$T^*M\oplus \R$. For those familiar with \cite{marius2} (to which we refer also
for notations), here are the main modifications: the splitting to be used is $(\sigma, 0)$,
the curvature becomes $(\Omega, \omega_r)$, the monodromy group is
$(A_{a}^{'}(r)\mathbb{Z}n, A_{a}(r)\mathbb{Z})$, and the function $r_{\mathcal{N}}$
will be $A_{a}{'}+ A_{a}$ away from the origin, and $+\infty$ at the origin. The conclusion
is that $M_a$ is integrable as a Jacobi manifold if and only if $\lim_{r\to 0}A'_a(r) +A_a(r)\neq 0$.
It is now easy to find various non-integrable examples:
\begin{enumerate}
\item[(i)] For $a(r)= r e^{r}$, $M_a$ is Poisson integrable but it is not Jacobi integrable.
\item[(ii)] For $a(r)=1/(\sin r +2)$, $M_a$ is Jacobi integrable but it is not Poisson integrable.
\item[(iii)] For $a$ such that $a(r)= r e^{\phi(r)r}$ where $\phi$ is a smooth function that equals to $1$ near the origin and equals
to $0$ for $r$ large enough, $M_a$ is neither Poisson nor Jacobi integrable.
\end{enumerate}
\end{ex}

\end{document}